\documentclass[11pt,a4paper,twoside,reqno,dvips]{article}
\usepackage{amssymb}
\usepackage [subnum]{cases}

\usepackage[top=30mm,bottom=30mm,left=31mm,right=31mm]{geometry}
\usepackage{comment}

\makeatletter 
\@addtoreset{equation}{section}
\makeatother  

\newtheorem{thm}{Theorem}[section]
\newtheorem{defn}[thm]{Definition}
\newtheorem{lem}[thm]{Lemma}
\newtheorem{cor}[thm]{Corollary}
\newtheorem{rem}[thm]{Remark}

\newtheorem{exam}[thm]{Example}

\begin{document}

\title{Stability of analytical and numerical  solutions of nonlinear  stochastic
delay differential equations\thanks{This work was partially supported by NSF of China (No.11171352, 11271311) and
       State Key Laboratory of High Performance Complex Manufacturing.}}


\author{{\bf Siqing Gan\thanks{Corresponding author. E-mail address: sqgan@csu.edu.cn}}\\
School of Mathematics and Statistics\\
State Key Laboratory of High Performance Complex Manufacturing\\
Central South University, Changsha, Hunan 410083, China\\
 \vspace{0.06cm}\\
 {\bf Aiguo Xiao}\\
 School of Mathematics and Computational Science\\
 Hunan Key Laboratory for Computation
      and Simulation in Science and Engineering\\
      Xiangtan University, Xiangtan, Hunan 411105, China\\
\vspace{0.06cm}\\
       {\bf Desheng Wang}\\
 Division of Mathematical Sciences, School of Physical and Mathematical Sciences\\
 Nanyang Technological University, Singapore, 637371}

\date{}

\maketitle

\begin{abstract}
\noindent This paper concerns the stability of analytical and numerical
solutions of nonlinear stochastic delay differential
equations (SDDEs). We derive sufficient conditions for the stability, contractivity and asymptotic contractivity in mean square of
the solutions for nonlinear SDDEs. The results provide a unified
theoretical treatment for SDDEs with constant delay and variable delay (including bounded and unbounded variable delays). Then the stability,
contractivity and asymptotic contractivity in mean square are investigated for
 the backward Euler method. It is shown that the backward Euler method preserves the properties of the underlying SDDEs.  The main results obtained in this work are different
 from those of Razumikhin-type theorems. Indeed,
our results hold without the necessity of constructing of finding an appropriate Lyapunov functional.

\vspace{0.1cm}

\vspace{0.1cm}

\end{abstract}

\section{Introduction}
\label{introduction}

Many physical, engineering and economic processes
can be modeled by stochastic differential equations (SDEs). The
rate of change of such a system depends only on its present state
and some noisy input. However, in many practical situations the rate of
change of the state depends not only on the present but also on
the past states of the system. Stochastic functional differential
equations (SFDEs) give a mathematical formulation for such system.
For more details on SFDEs, we refer to
\cite{ Mao97,Mohammed1984, Mohammed1996} and the references
therein.

SFDEs also can be regarded as a generalization of deterministic
functional differential equations when stochastic effects are taken
into account. For deterministic Volterra functional differential
equations (VFDEs) in Banach
spaces, Li \cite{Lishoufu2005} discussed the stability,
contractivity and asymptotic stability of the solutions. In \cite{Lishoufu2005}, the author introduced a so-called $\frac{1}{n}$-perturbed problem and constructed an auxiliary function $Q(t)$ for the corresponding study. The $\frac{1}{n}$-perturbed problem can be used to deal with a wide variety of delay arguments and the auxiliary function $Q(t)$ is the crux to establish the main results in \cite{Lishoufu2005}.  The work  \cite{Lishoufu2005} provides a unified framework for
stability analysis of nonlinear stiff problems in
ordinary differential equations, delay differential
equations, integro-differential equations and VFDEs
of other types. The theory in \cite{Lishoufu2005} was further extended to
nonlinear Volterra neutral functional differential equations
(VNFDEs) \cite{Wang2008}. Moreover, in \cite{WZ2010}, it is proved that the implicit Euler
method preserves the stability of VFDEs and VNFDEs.

It is natural to ask whether the solutions of SFDEs possess similar properties  to those presented in \cite{Lishoufu2005} and which methods can reproduce the properties. Due to the unique features of stochastic calculus, the
numerical analyses of SFDEs significantly differ from those developed for the numerical analyses of their deterministic
counterparts. In the literature, much attention on numerical stability has been focused
on a special class of SFDEs, namely, stochastic delay
differential equations (SDDEs); see \cite{Baker-Buckwar05,LCF04,Wang-Chen-11,Wang-Zhang-06,WMS2010,ZGH09}. The results mainly concern the
mean-square stability, asymptotic stability and
exponential mean-square stability for SDDEs with bounded lags.  Very recently, Fan, Song
and Liu \cite{Fan-Song-Liu-09} discussed the mean-square stability
of semi-implicit Euler methods for linear stochastic pantograph
equations. Far less is known for long-run
behavior of nonlinear SDDEs with unbounded lags. Moreover, to our best knowledge, there is no
work on the contractivity analysis of numerical methods for SDDEs.
Our aims in this paper are to investigate the
stability and contractivity of nonlinear SDDEs with bounded and
unbounded lags and to study the numerical preservation of those
the properties. The main results of this paper could be summarized as
follows.
\begin{enumerate}
\item[(i)] Sufficient conditions for the stability,
contractivity and asymptotic contractivity in mean square of the
solutions for nonlinear SDDEs are derived. The results provide a unified
theoretical treatment for SDDEs with constant delay and variable delay (including bounded and
unbounded variable delays). Applicability of the theory is illustrated by linear and nonlinear
SDDEs with a wide variety of delay arguments such as constant delays,
piecewise constant arguments, proportional delays and so on.  The theorems established in this
paper work for some SDDEs to which the existing theories cannot be
applied. Our main results of analytic solutions can
be regarded as a generalization of those in \cite{Lishoufu2005}
restricted in finite-dimensional Hilbert spaces and finitely many
delays to the stochastic version.

\item[(ii)] It is proved that the backward Euler method preserves the stability, contractivity and asymptotic contractivity in mean square of the underlying systems. In particular, Theorem \ref{thm-add-1} and Theorem \ref{thm-add-2} show that
the backward Euler method preserves the contractivity and asymptotic contractivity without any constraint on the numerical stepsize.

\end{enumerate}

We point out that the main theorems in the present paper are
different from the Razumikhin-type theorems established in \cite{Baker-Buckwar05,Mao97}.
Our theorems can be directly applied to establish the stability without the necessity of constructing and finding an appropriate Lyapunov functional, as required by the Razumikhin-type theorems.
In this sense, our theorems are more convenient for stability analysis than the Razumikhin-type theorems.

The rest of the paper is organized as follows. In section 2, we
introduce some notations and assumptions, which
will be used throughout the rest of the paper. In section 3,
some criteria for the stability, contractivity and asymptotic contractivity in mean square of solutions for nonlinear
SDDEs are established. The main results obtained in this section
are applied to SDDEs with bounded and unbounded lags,
respectively. In section 4, sufficient conditions for the stability, contractivity and asymptotic contractivity in mean square for the backward Euler method are derived. Stability of analytical and
numerical solutions of SDDEs with several delays is discussed in
section 5.

\section{Stochastic delay differential equations}

Let $(\Omega,\mathcal {F},\{\mathcal {F}_t\}_{t\geq
a},\mathbb{P})$ be a complete probability space with a filtration
$\{\mathcal {F}_t\}_{t\geq a}$ satisfying the usual conditions
(i.e., it is right continuous and $\mathcal {F}_{a}$ contains all
the $\mathbb{P}$-null sets). Let $w(t)=(w_1(t),...,w_m(t))^{T}$
be an $m$-dimensional Wiener process defined on the probability
space. Let $\langle \cdot\rangle$ be inner product in $ \mathbb{C}^d $  and $|\cdot|$
corresponding norm. In this paper, $|\cdot|$ also denotes the trace norm (F-norm) in
$\mathbb{C}^{d\times m}$. Also, $C([t_1,t_2];\mathbb{C}^d)$ is
used to represent the family of continuous mappings $\psi$ from
$[t_1,t_2]$ to $ \mathbb{C}^d $. Let $p>2$ and denote by $L^p_{\mathcal {F}_t}
([t_1,t_2];\mathbb{C}^d)$ the family of
${\mathcal {F}_t}$ -measurable $C([t_1,t_2];
\mathbb{C}^d)$-valued random variables $\psi=\{\psi(u): t_1 \leq
u \leq t_2\}$ such that $\|\psi\|^p_{\mathbb{E}} = \sup_{t_1 \leq u \leq
t_2}\mathbb{E}|\psi(u)|^p <\infty$. $\mathbb{E}$ denotes mathematical expectation with respect to $\mathbb{P}$.

Consider the following initial value problems of
SDDEs in the sense of It\^{o}
\begin{numcases} {\label{vsfde1}}
dx(t)=f(t,x(t),x(t-\tau(t)))dt + g(t,x(t),x(t-\tau(t)))dw(t),  \ t \in[a,b],\qquad\label{vsfde1a}\\
x(t)=\xi(t), \ \  t\in [a-\tau_0,a], \ \xi\in L^p_{\mathcal {F}_a} ([a-\tau_0, a];\mathbb{C}^d),\label{vsfde1b}
\end{numcases}
where  $a, b, \tau_0$ are constants with $-\infty<a<b<+\infty$ and $\tau_0 \geq
0, \tau(t)\geq 0, \inf\limits_{a\leq t\leq b}(t-\tau(t))\geq
a-\tau_0$, $f:[a,b]\times \mathbb{C}^d\times
\mathbb{C}^d\rightarrow\mathbb{C}^d, g:[a,b]\times
\mathbb{C}^d\times \mathbb{C}^d\rightarrow\mathbb{C}^{d\times m}$
are given continuous mappings.
We assume that the drift coefficient $f$ and the diffusion coefficient $g$ satisfy the following
conditions.
\begin{eqnarray}
&&\mbox{For each  $R>0$ there exists a constant
$C_R$,  depending only on  $R$,  such that
}\nonumber\\[-1.5ex] \label{vsfde-add-21}\\[-1.5ex]
 &&\qquad |f(t,x_1,y)-f(t,x_2,y)| \leq  C_R|x_1-x_2|
,\ \ \ |x_1|\vee
|x_2|\vee |y|\leq R,\nonumber\\
&&  \Re\big\langle x_1-x_2,
f(t,x_1,y)-f(t,x_2,y)\big\rangle \leq \alpha(t)|x_1-x_2|^2,\label{vsfde3}\\
&&  |f(t,x,y_1)-f(t,x,y_2)|\leq \beta(t)|y_1-y_2|, \label{vsfde4}\\
&&  |g(t,x_1,y_1)-g(t,x_2,y_2)|\leq
\gamma_1(t)|x_1-x_2|+\gamma_2(t)|y_1-y_2|, \label{vsfde5}
\end{eqnarray}
for all $t\in [a,b], x,x_1,x_2,y,y_1,y_2\in \mathbb{C}^d$, where $\Re a$ denotes the real part of the complex number $a$. Here $\alpha(t), \beta(t), \gamma_1(t)$ and $ \gamma_2(t)$ are
continuous real-valued functions. We introduce the following notations
\begin{eqnarray*}
&\mu_1^{(0)}=\inf\limits_{a\leq t\leq b}\tau(t)\geq 0,\ \ &
\mu_2^{(0)}(t_1,t_2)=\inf\limits_{t_1\leq t\leq
t_2}\big(t-\tau(t)\big)\geq a-\tau_0,\\
&&\qquad \forall t_1,t_2: a\leq
t_1\leq t_2\leq b.
\end{eqnarray*}
For convenience, we denote by
$\mathcal{SD}(\alpha,\beta,\gamma_1,\gamma_2)$ the all problems
(\ref{vsfde1}) which satisfy the conditions
(\ref{vsfde-add-21})-(\ref{vsfde5}). Such problems will be introduced in the next section (see Example \ref{example1} and Example \ref{example2}).

In order to deal with a wide variety of delay arguments, we introduce the so-called $\frac{1}{n}$-perturbed problem of (\ref{vsfde1}), which was first introduced by Li \cite{Lishoufu2005} for VFDEs.
We call the initial value problem
\begin{numcases}{\label{vsfde10}}
dx(t)=f\big(t,x(t),x^{(n,t)}(t-\tau(t))\big)dt + g\big(t,x(t),x^{(n,t)}(t-\tau(t))\big)dw(t),  t \in[a,b],\qquad\label{vsfde10a}\\
x(t)=\xi(t), \ \ t\in [a-\tau_0,a], \ \xi\in L^p_{\mathcal {F}_a} ([a-\tau_0, a];\mathbb{C}^d),  \label{vsfde10b}
\end{numcases}
an $\frac{1}{n}$-perturbed problem of the problem
(\ref{vsfde1}), where
\begin{equation}\label{nperturb1}
x^{(n,t)}(t-\tau(t))=\left\{
\begin{array}{ll}
x(t-\tau(t)),& \tau(t)\geq \frac{1}{n},\\
x(t-\frac{1}{n}),&\tau(t)<\frac{1}{n}.
\end{array}\right.
\end{equation}
Here the natural number $n>\frac{1}{\tau_0}$ can be
arbitrarily given. Without
lose of generality,  we always assume $\tau_0>0$. In fact, in the case of $\tau_0=0$, we can
replace $\tau_0$ by some positive number $\tilde{\tau}_0$ and define
$\xi(u)=\xi(a)$ for $u \in [a-\tilde{\tau}_0,a]$.

It is easy to verify that, if problem (\ref{vsfde1}) $\in
\mathcal{SD}(\alpha,\beta,\gamma_1,\gamma_2)$, then its
$\frac{1}{n}$-perturbed problem (\ref{vsfde10}) $\in
\mathcal{SD}(\alpha,\beta,\gamma_1,\gamma_2)$. It is clear that
$\tilde{\tau}(t)=\max\big\{\tau(t),\frac{1}{n}\big\}$ is the time-lag argument of the $\frac{1}{n}$-perturbed problem and
\begin{equation}\label{vsfde12}
\tilde{\mu}_1^{(0)}=\inf\limits_{a\leq t\leq b}\tilde{\tau}(t)\geq
\frac{1}{n}.
\end{equation}

\section{Stability analysis of SDDEs}

We discuss the following types of stability of SDDEs.

\begin{defn}
The solution of problem (\ref{vsfde1}) is said to be stable in mean square if
\begin{equation}\label{eqn-add-10}
\mathbb{E}|x(t)-y(t)|^2\leq C\sup\limits_{a-\tau_0\leq \theta\leq
a}\mathbb{E}|\xi(\theta)-\eta(\theta)|^2,
\end{equation}
where $y(t)$ is the solution of the perturbed problem
\begin{numcases} {\label{vsfde2-2}}
dy(t)=f(t,y(t),y(t-\tau(t)))dt + g(t,y(t),y(t-\tau(t)))dw(t),  \ t \in[a,b],\quad\label{vsfde2-2a}\\
y(t)=\eta(t),\label{vsfde2-2b} \, t\in [a-\tau_0,a], \ \ \eta\in
L^p_{\mathcal {F}_a} ([a-\tau_0, a];\mathbb{C}^d).
\end{numcases}
\end{defn}

\begin{defn}
The solution of problem (\ref{vsfde1}) is said to be
contractive in mean square if (\ref{eqn-add-10}) with $C\leq 1$ holds.
\end{defn}

\begin{defn}
The solution of problem (\ref{vsfde1}) is said to be asymptotically contractive in mean square if
$$
\lim\limits_{t\rightarrow +\infty}\mathbb{E}|x(t)-y(t)|^2=0,
$$
for which $[a,b]$ is replaced by $[a,+\infty)$ in (\ref{vsfde1a}) and (\ref{vsfde2-2a}).
\end{defn}

\begin{rem}
In the strict sense, (\ref{eqn-add-10}) with $C\leq 1$ means generalized contraction. For brevity, we simply call the solution is contractive in mean square.

There exist well-known stability definitions in literatures which are closely related to those presented in this paper, but there are differences among them. The existing notions of stability include mean-square stability for SDEs, that is, $\lim\limits_{t\rightarrow +\infty}\mathbb{E}|x(t)|^2=0$ (cf. \cite{SM96}); exponential mean-square contraction of trajectories for SDEs with jumps (cf. \cite{Higham-Kloeden-2005}). The contractivity in mean square is weaker than that in \cite{Higham-Kloeden-2005}.
\end{rem}

The continuity of $f$ and $g$ implies that
\begin{equation}\label{vsfde-add-5}
|f(t,0,0)|\leq C, \ \ |g(t,0,0)|\leq C,\ \ t\in [a,b],
\end{equation}
where $C$ only depends on $f,g$ and the  interval $[a,b]$.
We note that condition (\ref{vsfde5}) implies that the diffusion coefficient $g$ satisfies the local Lipschitz condition
\begin{eqnarray}
& |g(t,x_1,y_1)-g(t,x_2,y_2)|
 \leq  C_R\big(|x_1-x_2|
+|y_1-y_2|\big),\nonumber\\[-1.5ex]\label{vsfde-add-54}\\[-1.5ex]
&\qquad \qquad t\in [a,b], x_1,x_2,y_1,y_2\in \mathbb{C}^d,  \
|x_1|\vee |x_2|\vee |y_1|\vee |y_2|\leq R.\nonumber
\end{eqnarray}
In fact, we can choose any $C_R$ with $C_R\geq \max\{\max\limits_{a\leq t\leq b}\gamma_1(t), \max\limits_{a\leq t\leq b}\gamma_2(t)\}$.
Using (\ref{vsfde3}), (\ref{vsfde4}) and (\ref{vsfde-add-5}), we
have
\begin{eqnarray}
\Re \langle x, f(t,x,y) \rangle &=&\Re \Big\langle x-0,
f(t,x,y)-f(t,0,y)+f(t,0,y)-f(t,0,0)+f(t,0,0)
\Big\rangle\nonumber\\
&\leq & \alpha (t)|x|^2+\beta(t)|x||y|+C|x|
\leq  C_1\Big(1+|x|^2+|y|^2\Big),\label{vsfde-add-24}
\end{eqnarray}
where $C_1$ only depends on $C,\ \max\limits_{a\leq t\leq b}\alpha(t)$
and $\max\limits_{a\leq t\leq b}\beta(t)$.
By (\ref{vsfde5}), we have
\begin{equation}\label{vsfde-add-25}
|g(t,x,y)|^2 \leq  2|g(t,x,y)-g(t,0,0)|^2+2|g(t,0,0)|^2
 \leq C_1\Big(1+|x|^2+|y|^2\Big),
\end{equation}
where $C_1$ only depends on $C, \ \max\limits_{a\leq t\leq
b}\gamma_1^2(t)$ and $\max\limits_{a\leq t\leq b}\gamma_2^2(t)$.

\subsection{Finite interval}

In order to prove the main theorems in this section, we prepare
the following lemmas.

\begin{lem}\label{lem-add-4}
Assume that problem (\ref{vsfde1}) $\in
\mathcal{SD}(\alpha,\beta,\gamma_1,\gamma_2)$. Then for each
$p\geq 2$ there is
$\bar{C}=\bar{C}(p,a,b,\alpha,\beta,\gamma_1,\gamma_2)$ such that
\begin{equation}\label{vsfde-add-26}
\mathbb{E}\big(\sup\limits_{a\leq t\leq b}|x(t)|^p\big)\leq
\bar{C}\big(1+\mathbb{E} (\sup\limits_{a-\tau_0\leq t\leq
a}|\xi(t)|^p)\big)=A.
\end{equation}
\end{lem}

\noindent{\bf Proof.}
For every integer $k\geq 1$, define the stopping
time
\begin{equation}
\rho_k= \inf\{t\in [a,b]: \sup\limits_{a-\tau_0\leq
\theta\leq t}|x(\theta)|\geq k\},
\end{equation}
where we use the convention $\rho_k=b$ if the set is empty in the right-hand side. Clearly,
$\rho_k\uparrow b$ almost surely as $k\rightarrow +\infty$. Let
$x^k(t)=x(t\wedge \rho_k)$.
Using the It\^{o} formula, we have for $a\leq t\leq b$
\begin{eqnarray*}
&&\Big(1+|x^k(t)|^2\Big)^{\frac{p}{2}}
=\Big(1+|\xi(a)|^2\Big)^{\frac{p}{2}}\\
&&+p\int_a^t
\Big(1+|x^k(s)|^2\Big)^{\frac{p-2}{2}}\Re \langle x^k(s),
f(s,x^k(s),x^k(s-\tau(s))) \rangle I_{[[a,\rho_k]]}(s)ds   \\
&&+\frac{p}{2}\int_a^t \Big(1+|x^k(s)|^2\Big)^{\frac{p-2}{2}}
|g(s,x^k(s),x^k(s-\tau(s)))|^2 I_{[[a,\rho_k]]}(s)ds\\
&&+\frac{p(p-2)}{2}\int_a^t
\Big(1+|x^k(s)|^2\Big)^{\frac{p-4}{2}}|(x^k(s))^Tg(s,x^k(s),x^k(s-\tau(s)))|^2I_{[[a,\rho_k]]}(s)ds   \\
\qquad &&+p\int_a^t \Big(1+|x^k(s)|^2\Big)^{\frac{p-2}{2}}\Re \langle
x^k(s), g(s,x^k(s),x^k(s-\tau(s))) \rangle I_{[[a,\rho_k]]}(s)dw(s),
\end{eqnarray*}
by (\ref{vsfde-add-24}) and
(\ref{vsfde-add-25}), and hence
\begin{eqnarray*}
&&\Big(1+|x^k(t)|^2\Big)^{\frac{p}{2}}\\
&\leq &
2^{\frac{p-2}{2}}\Big(1+|\xi(a)|^p\Big)+2pC_1\int_a^t
\Big(1+|x^k(s)|^2\Big)^{\frac{p-2}{2}}\Big(1+\sup\limits_{s-\tau(s)\leq
u\leq s}|x^k(u)|^2\Big)ds\\
&&+pC_1\int_a^t
\Big(1+|x^k(s)|^2\Big)^{\frac{p-2}{2}}\Big(1+\sup\limits_{s-\tau(s)\leq
u\leq s}|x^k(u)|^2\Big)ds\\
&&+p(p-2)C_1\int_a^t
\Big(1+|x^k(s)|^2\Big)^{\frac{p-4}{2}}|x^k(s)|^2\Big(1+\sup\limits_{s-\tau(s)\leq
u\leq s}|x^k(u)|^2\Big)ds\\
&&+p\int_a^t \Big(1+|x^k(s)|^2\Big)^{\frac{p-2}{2}}\Re \langle
x^k(s), g(s,x^k(s),x^k(s-\tau(s))) \rangle I_{[[a,\rho_k]]}(s)dw(s)\\
&\leq & 2^{\frac{p-2}{2}}\Big(1+|\xi(a)|^p\Big)
+p(p+1)C_1\int_a^t\Big(1+\sup\limits_{s-\tau(s)\leq u\leq
s}|x^k(u)|^2\Big)^{\frac{p}{2}}ds\\
&&+p\int_a^t \Big(1+|x^k(s)|^2\Big)^{\frac{p-2}{2}}\Re \langle
x^k(s), g(s,x^k(s),x^k(s-\tau(s))) \rangle I_{[[a,\rho_k]]}(s)dw(s),
\end{eqnarray*}
which yields
\begin{eqnarray*}
&&\mathbb{E}\sup\limits_{a\leq s\leq
t}\big(1+|x^k(s)|^2\big)^{\frac{p}{2}}\\
& \leq &
2^{\frac{p-2}{2}}\Big(1+\mathbb{E}\sup\limits_{a-\tau_0\leq s\leq
a}|\xi(s)|^p\Big)+C_2\mathbb{E}\int_a^t\Big(1+\sup\limits_{a-\tau_0\leq
u\leq
s}|x^k(u)|^2\Big)^{\frac{p}{2}}ds  \\
&&+p\mathbb{E}\Big(\sup\limits_{a\leq s\leq t}\int_a^s
\big(1+|x^k(u)|^2\big)^{\frac{p-2}{2}}\Re \langle x^k(u),
g(u,x^k(u),x^k(u-\tau(u))) \rangle I_{[[a,\rho_k]]}(u)dw(u)\Big).
\end{eqnarray*}
Applying the Burkholder-Davis-Gundy inequality to the third term
on the right-hand side of the above inequality, we obtain the bound
\begin{eqnarray*}
&&p\mathbb{E}\Big(\sup\limits_{a\leq s\leq t}\int_a^s
\big(1+|x^k(u)|^2\big)^{\frac{p-2}{2}}\Re \langle x^k(u),
g(u,x^k(u),x^k(u-\tau(u))) \rangle I_{[[a,\rho_k]]}(u)dw(u)\Big)  \\
&\leq &
C_p\mathbb{E}\Big(\int_a^t\big(1+|x^k(u)|^2\big)^{p-2}|x^k(u)|^2|g(u,x^k(u),x^k(u-\tau(u)))|^2du
\Big)^{\frac{1}{2}}  \\
&\leq & C_p\mathbb{E}\Big(\sup\limits_{a\leq u\leq
t}\big(1+|x^k(u)|^2\big)^{\frac{p}{2}}\int_a^t\big(1+|x^k(u)|^2\big)^{\frac{p-4}{2}}
|x^k(u)|^2|g(u,x^k(u),x^k(u-\tau(u)))|^2du \Big)^{\frac{1}{2}}  \\
&\leq & \frac{1}{2}\mathbb{E}\sup\limits_{a\leq u\leq
t}\big(1+|x^k(u)|^2\big)^{\frac{p}{2}}+\frac{C_p^2}{2}2C_1\mathbb{E}
\int_a^t\Big(1+\sup\limits_{a-\tau_0\leq s\leq
u}|x^k(s)|^2\Big)^{\frac{p}{2}}du.
\end{eqnarray*}
Consequently,
\begin{eqnarray*}
\mathbb{E}\sup\limits_{a\leq s\leq
t}\big(1+|x^k(s)|^2\big)^{\frac{p}{2}} &\leq&
2^{\frac{p}{2}}\big(1+\mathbb{E}\sup\limits_{a-\tau_0\leq s\leq
a}|\xi(s)|^p\big)+C_3\mathbb{E}\int_a^t\Big(1+\sup\limits_{a-\tau_0\leq
s\leq u}|x^k(s)|^2\Big)^{\frac{p}{2}}du  \\
&=&2^{\frac{p}{2}}\big(1+\mathbb{E}\sup\limits_{a-\tau_0\leq s\leq
a}|\xi(s)|^p\big)+C_3\mathbb{E}\int_a^t\sup\limits_{a-\tau_0\leq
s\leq u}\Big(1+|x^k(s)|^2\Big)^{\frac{p}{2}}du.
\end{eqnarray*}
Further, we notice that
\begin{eqnarray*}
\mathbb{E}\sup\limits_{a-\tau_0\leq s\leq
t}\big(1+|x^k(s)|^2\big)^{\frac{p}{2}} &\leq &
\mathbb{E}\sup\limits_{a-\tau_0\leq s\leq
a}\big(1+|\xi(s)|^2\big)^{\frac{p}{2}}+\mathbb{E}\sup\limits_{a\leq
s\leq t}\big(1+|x^k(s)|^2\big)^{\frac{p}{2}}  \\
&\leq&  2^{\frac{p-2}{2}}(1+\mathbb{E}\sup\limits_{a-\tau_0\leq s\leq
a}|\xi(s)|^p)+\mathbb{E}\sup\limits_{a\leq s\leq
t}\big(1+|x^k(s)|^2\big)^{\frac{p}{2}}.
\end{eqnarray*}
Therefore,
\begin{eqnarray*}
\mathbb{E}\sup\limits_{a-\tau_0\leq s\leq
t}\big(1+|x^k(s)|^2\big)^{\frac{p}{2}}\leq
\frac{3}{2}2^{\frac{p}{2}}(1+\mathbb{E}\sup\limits_{a-\tau_0\leq
s\leq a}|\xi(s)|^p)+C_3\int_a^t\mathbb{E}\sup\limits_{a-\tau_0\leq
s\leq u}\Big(1+|x^k(s)|^2\Big)^{\frac{p}{2}}du.
\end{eqnarray*}
Now the Gronwall's inequality yields that
\begin{eqnarray*}
\mathbb{E}\big(\sup\limits_{a-\tau_0\leq s\leq t}|x^k(s)|^p\big)\leq
\mathbb{E}\big(\sup\limits_{a-\tau_0\leq s\leq
t}(1+|x^k(s)|^2)^{\frac{p}{2}}\big)\leq \frac{3}{2}2^{\frac{p}{2}}
\big(1+\mathbb{E}\sup\limits_{a-\tau_0\leq s\leq
a}|\xi(s)|^p\big)e^{C_3(t-a)}.
\end{eqnarray*}
Letting $k\rightarrow +\infty$ and applying the Fatou's lemma, we
obtain the desired result.

\begin{lem}\label{lem-add-3}
Assume that problem (\ref{vsfde1}) $\in
\mathcal{SD}(\alpha,\beta,\gamma_1,\gamma_2)$. Then there exists a
unique solution $x(t)$ to equation (\ref{vsfde1}).
\end{lem}

\noindent{\bf Proof.}
Using (\ref{vsfde-add-21}), (\ref{vsfde4}), (\ref{vsfde-add-54}),
(\ref{vsfde-add-24}) and (\ref{vsfde-add-25}), we are able to complete the
proof along the lines of the ones for Theorem 5.2.7 and Theorem
2.3.5 in~\cite{Mao97}.

\begin{lem}\label{lem-add-2}
Assume that problem (\ref{vsfde1}) $\in
\mathcal{SD}(\alpha,\beta,\gamma_1,\gamma_2)$. Then we have
\begin{equation}\label{vsfde-add-35}
\lim\limits_{v\rightarrow 0}\sup\limits_{|t-s|=v}\mathbb{E}|x(t)-x(s)|^2=0,\ \ a\leq s,t\leq b.
\end{equation}
\end{lem}

\noindent{\bf Proof.} The proof is mainly based on the techniques employed in the
proof of Theorem 2.2 in \cite{HMS2002}.
Integrating (\ref{vsfde1}) gives for $a\leq s<t\leq b$
\begin{equation}\label{vsfde-add-27}
x(t)-x(s)=\int_s^t f(u,x(u),x(u-\tau(u)))du+\int_s^t
g(u,x(u),x(u-\tau(u)))dw(u).
\end{equation}
Let $e(s,t)=x(t)-x(s)$,
$$
\rho_R=\inf\{t\in [a,b]:
\sup\limits_{a-\tau_0\leq \theta\leq t}|x(\theta)|\geq R\},
$$
where $\rho_R=b$ if the set is empty in the right-hand side.
Using
the Young inequality: for $r^{-1}+q^{-1}=1$
$$
ab\leq \frac{\delta}{r}a^r+\frac{1}{q\delta^{q/r}}b^q, \ \forall a,b,\delta>0
$$
and letting $r=\frac{p}{2}, q=\frac{p}{p-2}$, we thus have for any $\delta>0$
\begin{eqnarray}
&&\mathbb{E}\big(|e(s,t)|^2\big)=\mathbb{E}\big(|e(s,t)|^2I_{\{\rho_R\geq b\}}\big)+\mathbb{E}\big(|e(s,t)|^2I_{\{\rho_R<b\}}\big)\nonumber
\\[-1.5ex]
\label{vsfde-add-29}\\[-1.5ex]
&\leq&
\mathbb{E}\big(|e(s,t)|^2I_{\{\rho_R\geq b\}}\big)+\frac{2\delta}{p}\mathbb{E}\big(|e(s,t)|^p\big)
+\frac{1-\frac{2}{p}}{\delta^{2/(p-2)}}\mathbb{P}(\rho_R<b),\nonumber
\end{eqnarray}
where $p>2$.
It follows from (\ref{vsfde-add-26}) that
\begin{eqnarray}
\mathbb{P}(\rho_R<b)=\mathbb{E}\left(I_{\{\rho_R<b\}}\frac{|x(\rho_R)|^p}{R^p}\right)\leq
\frac{1}{R^p}\mathbb{E}\big(\sup\limits_{a\leq t\leq
b}|x(t)|^p\big)\leq \frac{A}{R^p},\label{vsfde-add-30}\\
\mathbb{E}\big(|e(s,t)|^p\big)\leq
2^{p-1}\mathbb{E}\big(\sup\limits_{a\leq s\leq
b}|x(s)|^p+\sup\limits_{a\leq t\leq b}|x(t)|^p)\big)\leq 2^pA.\label{vsfde-add-31}
\end{eqnarray}
We then have
\begin{equation}\label{vsfde-add-32}
\mathbb{E}\big(|e(s,t)|^2\big)\leq
\mathbb{E}\big(|e(s,t)|^2I_{\{\rho_R\geq b\}}\big)+
\frac{2^{p+1}\delta A}{p}+\frac{(p-2)A}{p\delta^{2/(p-2)}R^p}.
\end{equation}
Further, using the H$\ddot{\mbox{o}}$lder's inequality and the It\^{o}
isometry, we obtain
\begin{eqnarray*}
&&\mathbb{E}\big(|e(s,t)|^2I_{\{\rho_R\geq b\}}\big)\\
&=&\mathbb{E}\Big(\Big|\int_s^tf(u,x(u),x(u-\tau(u)))du+\int_s^tg(u,x(u),x(u-\tau(u)))dw(u)\Big|^2I_{\{\rho_R\geq b\}}\Big)\\
&\leq
&2\mathbb{E}\Big(\Big(\big|\int_s^tf(u,x(u),x(u-\tau(u)))du\big|^2+\big|\int_s^tg(u,x(u),x(u-\tau(u)))dw(u)\big|^2\Big)
I_{\{\rho_R\geq b\}}\Big)\\
&\leq & 2(t-s)\mathbb{E}\Big(I_{\{\rho_R\geq b\}}\int_s^t|f(u,x(u),x(u-\tau(u)))|^2du\Big)+2\mathbb{E}\int_s^t|g(u,x(u),x(u-\tau(u)))|^2du\\
&\leq & 2(t-s)\mathbb{E}\Big(I_{\{\rho_R\geq b\}}\int_s^t|f(u,x(u),x(u-\tau(u)))-f(u,0,0)+f(u,0,0)|^2du\Big)\\
&&+2\mathbb{E}\int_s^t|g(u,x(u),x(u-\tau(u)))|^2du.\\
\end{eqnarray*}
By (\ref{vsfde-add-21}), (\ref{vsfde4}), (\ref{vsfde-add-5}), (\ref{vsfde-add-25})
and Lemma \ref{lem-add-4}, we have
\begin{eqnarray*}
&&\mathbb{E}\big(|e(s,t)|^2I_{\{\rho_R\geq b\}}\big)  \\
&\leq &
8C_R^2(t-s)\int_s^t\Big(\mathbb{E}|x(u)|^2+\mathbb{E}|x(u-\tau(u))|^2\Big)du
+4(t-s)\mathbb{E}\int_s^t |f(u,0,0)|^2du  \\
&&+2C_1\int_s^t\Big(1+\mathbb{E}|x(u)|^2+\mathbb{E}|x(u-\tau(u))|^2\Big)du
\leq  C_4(t-s),
\end{eqnarray*}
where $C_4$ is independent of $s$ and $t$.
A combination of this expression and (\ref{vsfde-add-32}) leads to
\begin{equation}\label{vsfde-add-34}
\mathbb{E}\big(|e(s,t)|^2\big)\leq C_4(t-s)+ \frac{2^{p+1}\delta
A}{p}+\frac{(p-2)A}{p\delta^{2/(p-2)}R^p}.
\end{equation}
Therefore, for any given $\epsilon>0$, we can choose $\delta$ and $R$ such that
\begin{equation}
\frac{2^{p+1}\delta A}{p}<\frac{1}{3}\epsilon,\ \ \frac{(p-2)A}{p\delta^{2/(p-2)}R^p}<\frac{1}{3}\epsilon,
\end{equation}
and then choose $t-s$ sufficiently small such that
$C_4(t-s)<\frac{1}{3}\epsilon$. Hence, we have
\begin{eqnarray*}
\lim\limits_{v\rightarrow 0}\sup\limits_{|t-s|=v}\mathbb{E}\big(|x(t)-x(s)|^2\big)=\lim\limits_{v\rightarrow 0}\sup\limits_{|t-s|=v}\mathbb{E}\big(|e(s,t)|^2\big)=0.
\end{eqnarray*}
The proof is complete.

\begin{lem}\label{lem1}
Assume that problem (\ref{vsfde1}) $\in
\mathcal{SD}(\alpha,\beta,\gamma_1,\gamma_2)$. Then for any $t_1,
t_2: a\leq t_1\leq t_2\leq b$,
\begin{equation}\label{vsfde2-1}
\begin{array}{l}
\mathbb{E}|x(t_2)-y(t_2)|^2\leq
e^{\int_{t_1}^{t_2}\sigma(t)dt}\mathbb{E}|x(t_1)-y(t_1)|^2\\
\qquad +\int_{t_1}^{t_2}\varrho(s)e^{\int_s^{t_2}
\sigma(u)du}ds\sup\limits_{\mu_2^{(0)}(t_1,t_2)\leq \theta\leq
t_2-\mu_1^{(0)}}\mathbb{E}|x(\theta)-y(\theta)|^2,
\end{array}
\end{equation}
where $y(t)$ is the solution of the perturbed problem (\ref{vsfde2-2}), and $\sigma(t), \varrho(t)$ are defined by
\begin{equation}\label{denote1}
\begin{array}{ll}
&\sigma(t)=2\alpha(t)+\beta(t)+\gamma_1(t)\gamma_2(t)+\gamma_1^2(t),\\
&\varrho(t)=\beta(t)+\gamma_1(t)\gamma_2(t)+\gamma_2^2(t).
\end{array}
\end{equation}
\end{lem}

\noindent{\bf Proof.}
Let
\begin{equation}\label{vsfde2-3}
V(t,x(t))=p(t)\big(|x(t)|^2+\delta q(t)\big),\ \ t_1\leq t\leq
t_2,
\end{equation}
where
\begin{equation}\label{vsfde2-4}
p(t)=e^{-\int_{a}^t\sigma(u)du},\ \ q(t)=-(p(t))^{-1}\int_{a}^t\varrho(u)p(u)du,
\end{equation}
$\delta$ is a constant to be determined.  Then we have
$ p'(t)=-\sigma(t)p(t),
\big(p(t)q(t)\big)'=-p(t)\varrho(t)$.
By  (\ref{vsfde2-3}),(\ref{vsfde2-4})
and the It\^{o} formula,
one can derive that, for $a\leq t_1\leq t_2\leq b$,
\begin{eqnarray*}
&&\mathbb{E}V\big(t_2,x(t_2)-y(t_2)\big)=\mathbb{E}V\big(t_1,x(t_1)-y(t_1)\big)\\
&&+\int_{t_1}^{t_2}\bigg\{-\sigma(t)p(t)\mathbb{E}|x(t)-y(t)|^2-\delta
p(t)\varrho(t)\bigg. \\
&& \bigg.+2p(t)\mathbb{E}\Re\big\langle
x(t)-y(t),f(t,x(t),x(t-\tau(t)))-f(t,y(t),y(t-\tau(t)))\big\rangle\\
&&+p(t)\mathbb{E}\big|g(t,x(t),x(t-\tau(t)))-g(t,y(t),y(t-\tau(t)))\big|^2\bigg\}dt\\
&&\leq \mathbb{E}V\big(t_1,x(t_1)-y(t_1)\big)+\int_{t_1}^{t_2}\bigg\{-\sigma(t)p(t)\mathbb{E}|x(t)-y(t)|^2-\delta
p(t)\varrho(t)\bigg.\\
&&+2p(t)\mathbb{E}\Re \big\langle x(t)-y(t),f(t,x(t),x(t-\tau(t)))-f(t,y(t),x(t-\tau(t)))\big\rangle\\
&&\bigg.+2p(t)\mathbb{E}\Re\big\langle x(t)-y(t),
f(t,y(t),x(t-\tau(t)))-f(t,y(t),y(t-\tau(t)))\big\rangle\\
&&+p(t)\mathbb{E}|g(t,x(t),x(t-\tau(t)))-g(t,y(t),y(t-\tau(t)))|^2\bigg\}dt.
\end{eqnarray*}
Using the conditions (\ref{vsfde3})-(\ref{vsfde5}), we obtain
\begin{eqnarray*}
&&\mathbb{E}V\big(t_2,x(t_2)-y(t_2)\big)\leq \mathbb{E}V\big(t_1,x(t_1)-y(t_1)\big)  \\
&&+\int_{t_1}^{t_2}
\bigg\{-\sigma(t)p(t)\mathbb{E}|x(t)-y(t)|^2-\delta
p(t)\varrho(t)  \bigg.\\
&&+2p(t)\alpha(t)\mathbb{E}|x(t)-y(t)|^2
+2p(t)\beta(t)\mathbb{E}\big(|x(t)-y(t)||x(t-\tau(t))-y(t-\tau(t))|\big)  \\
&&\bigg.+p(t)\mathbb{E}\Big(\gamma_1(t)|x(t)-y(t)|+\gamma_2(t)|x(t-\tau(t))-y(t-\tau(t))|\Big)^2\bigg\}dt  \\
&&\leq \mathbb{E}V\big(t_1,x(t_1)-y(t_1)\big)+\int_{t_1}^{t_2}\bigg\{-\sigma(t)p(t)\mathbb{E}|x(t)-y(t)|^2-\delta
p(t)\varrho(t)\bigg.  \\
&&+2p(t)\alpha(t)\mathbb{E}|x(t)-y(t)|^2
+p(t)\beta(t)\Big(\mathbb{E}|x(t)-y(t)|^2+\mathbb{E}|x(t-\tau(t))-y(t-\tau(t))|^2
\Big)  \\
&&\bigg.+p(t)\gamma_1^2(t)\mathbb{E}|x(t)-y(t)|^2+p(t)\gamma_2^2(t)
\mathbb{E}|x(t-\tau(t))-y(t-\tau(t))|^2  \\
&&+\bigg.p(t)\gamma_1(t)\gamma_2(t)\Big(\mathbb{E}|x(t)-y(t)|^2+\mathbb{E}|x(t-\tau(t))-y(t-\tau(t))|^2\Big)
\bigg\}dt\\
&&\leq \mathbb{E}V\big(t_1,x(t_1)-y(t_1)\big)-\delta\int_{t_1}^{t_2}p(t)\varrho(t)dt\\
&&+\int_{t_1}^{t_2}p(t)\Big(-\sigma(t)+2\alpha(t)+\beta(t)+\gamma_1^2(t)+\gamma_1(t)\gamma_2(t)\Big)\mathbb{E}|x(t)-y(t)|^2dt\\
&&+\int_{t_1}^{t_2}p(t)\Big(\beta(t)+\gamma_2^2(t)+\gamma_1(t)\gamma_2(t)\Big)\mathbb{E}|x(t-\tau(t))-y(t-\tau(t))|^2dt.
\end{eqnarray*}
The substitution of (\ref{denote1}) into this  gives
\begin{equation}\label{vsfde2-6}
\begin{array}{rcl}
&&\mathbb{E}V\big(t_2,x(t_2)-y(t_2)\big)\\
&&\leq
\mathbb{E}V\big(t_1,x(t_1)-y(t_1)\big)+\int_{t_1}^{t_2}p(t)\beta(t)\Big(\mathbb{E}|x(t-\tau(t))-y(t-\tau(t))|^2
-\delta\Big)dt  \\
&&+\int_{t_1}^{t_2}p(t)\gamma_2^2(t)\Big(\mathbb{E}|x(t-\tau(t))-y(t-\tau(t))|^2-\delta\Big)dt  \\
&&+\int_{t_1}^{t_2}p(t)\gamma_1(t)\gamma_2(t)\Big(\mathbb{E}|x(t-\tau(t))-y(t-\tau(t))|^2-\delta\Big)dt.
\end{array}
\end{equation}
Lemma \ref{lem-add-4} implies that
$\mathbb{E}\sup\limits_{a-\tau_0\leq t\leq b}|x(t)|^2< +\infty,
\mathbb{E}\sup\limits_{a-\tau_0\leq t\leq b}|y(t)|^2< +\infty$.
Consequently,
$\sup\limits_{a-\tau_0\leq t\leq b}\mathbb{E}|x(t)-y(t)|^2< +\infty.$
Let $\delta=\sup\limits_{\mu_2^{(0)}(t_1,t_2)\leq t\leq
t_2-\mu_1^{(0)}}\mathbb{E}|x(t)-y(t)|^2.$
It follows from (\ref{vsfde2-6}) that
\begin{equation}\label{vsfde2-6a}
\mathbb{E}V\big(t_2,x(t_2)-y(t_2)\big)\leq
\mathbb{E}V\big(t_1,x(t_1)-y(t_1)\big).
\end{equation}
The required estimate (\ref{vsfde2-1}) now follows from
(\ref{vsfde2-6a}).

\begin{lem}\label{lem2}
Assume that problem (\ref{vsfde1}) $\in
\mathcal{SD}(\alpha,\beta,\gamma_1,\gamma_2)$. Then
\begin{equation}\label{vsfde2-8}
\lim\limits_{n\rightarrow \infty}\sup\limits_{a-\tau_0\leq t\leq
b}\mathbb{E}|x(t)-x_n(t)|^2=0,
\end{equation}
where $x_n(t)$ is the solution of the $\frac{1}{n}$-perturbed
problem (\ref{vsfde10}).
\end{lem}

\noindent{\bf Proof.}
For any given natural number $n>\frac{1}{\tau_0}$, we can choose a
natural number $q$ sufficiently large such that
$
\mu=(b-a)/q<\frac{1}{n}.
$
Let
\begin{eqnarray*}
&&t_1=a+(i-1)\mu,\ \ t_2=a+i\mu,\ \ i=1,2,\cdots,q,\\
&&\alpha_0=\max\Big\{\max\limits_{a\leq t\leq
b}\sigma(t),0\Big\},\ \beta_0=\max\limits_{a\leq t\leq b}\varrho(t),\ \gamma_0=\max\Big\{\max\limits_{a\leq t\leq
b}\big(\sigma(t)+\varrho(t)\big),1\Big\},\\
&&V(t,x(t))=p(t)\big(|x(t)|^2+\delta q(t)\big), \ \ t_1\leq t\leq
t_2,
\end{eqnarray*}
where
\begin{equation}\label{vsfde2-10}
\delta=2(\varepsilon_n+\sup\limits_{a-\tau_0\leq \theta\leq
t_1}\mathbb{E}|x(\theta)-x_n(\theta)|^2),\ \
 \varepsilon_n=\sup\limits_{a\leq t\leq
b}\big(\sup\limits_{t-\frac{1}{n}\leq \theta\leq
t}\mathbb{E}\big|x(\theta)-x(t-\frac{1}{n})\big|^2\big),
\end{equation}
$p(t)$ and $q(t)$ are defined by (\ref{vsfde2-4}), $\sigma(t)$ and $\varrho(t)$ are defined by (\ref{denote1}).
For $t_1\leq t\leq t_2$, we can obtain the following estimate in the same way as (\ref{vsfde2-6})
\begin{equation}\label{vsfde-add-18}
\begin{array}{rcl}
&&\mathbb{E}V\big(t,x(t)-x_n(t)\big) \leq
\mathbb{E}V\big(t_1,x(t_1)-x_n(t_1)\big)  \\
&&+\int_{t_1}^tp(s)\beta(s)\Big(\mathbb{E} \big|x(s-\tau(s))-x_n^{(n,s)}(s-\tau(s))\big|^2-\delta\Big)ds  \\
&&+\int_{t_1}^tp(s)\gamma_2^2(s)\Big(\mathbb{E}
\big|x(s-\tau(s))-x_n^{(n,s)}(s-\tau(s))\big|^2-\delta\Big)ds  \\
&&+\int_{t_1}^tp(s)\gamma_1(s)\gamma_2(s)\Big(\mathbb{E}
\big|x(s-\tau(s))-x_n^{(n,s)}(s-\tau(s))\big|^2-\delta\Big)ds  \\
&&\leq
\mathbb{E}V\big(t_1,x(t_1)-x_n(t_1)\big)  \\
&&+\Big(\sup\limits_{a\leq u\leq
t}\mathbb{E}\big|x(u-\tau(u))-x_n^{(n,u)}(u-\tau(u))\big|^2-\delta\Big)\int_{t_1}^tp(s)
\varrho(s)ds.
\end{array}
\end{equation}
Moreover, by (\ref{nperturb1}) we find that
\begin{eqnarray*}
&&\sup\limits_{a\leq u\leq
t}\mathbb{E}\big|x(u-\tau(u))-x_n^{(n,u)}(u-\tau(u))\big|^2\\
&=& \max\Big\{\sup\limits_{a\leq u\leq
t\atop \tau(u)\geq \frac{1}{n}}\mathbb{E}|x(u-\tau(u))-x_n^{(n,u)}(u-\tau(u))|^2,\\
&&\qquad\quad \sup\limits_{a\leq u\leq
t\atop 0\leq \tau(u)<\frac{1}{n}}\mathbb{E}\big|x(u-\tau(u))-x_n^{(n,u)}(u-\tau(u))\big|^2\Big\}\\
&\leq & \max\Big\{\sup\limits_{a-\tau_0\leq u\leq
t-\frac{1}{n}}\mathbb{E}|x(u)-x_n(u)|^2,\sup\limits_{a\leq u\leq
t}\Big(\sup\limits_{u-\frac{1}{n}\leq s\leq u}\mathbb{E}\big|x(s)-x_n(u-\frac{1}{n})\big|^2\Big)\Big\}\\
&\leq & \max\Big\{\sup\limits_{a-\tau_0\leq u\leq
t-\frac{1}{n}}\mathbb{E}|x(u)-x_n(u)|^2,\\
&&\qquad \qquad\sup\limits_{a\leq u\leq
t}\Big(\sup\limits_{u-\frac{1}{n}\leq s\leq u}\mathbb{E}\big|x(s)-x(u-\frac{1}{n})+x(u-\frac{1}{n})-x_n(u-\frac{1}{n})\big|^2\Big)\Big\}\\
&\leq & 2\max\Big\{\sup\limits_{a-\tau_0\leq u\leq
t-\frac{1}{n}}\mathbb{E}|x(u)-x_n(u)|^2,\sup\limits_{a\leq u\leq
t}\Big(\sup\limits_{u-\frac{1}{n}\leq s\leq u}\mathbb{E}\big|x(s)-x(u-\frac{1}{n})\big|^2\Big)\Big\}\\
&\leq & 2\Big\{\sup\limits_{a-\tau_0\leq u\leq
t_1}\mathbb{E}|x(u)-x_n(u)|^2+\sup\limits_{a\leq u\leq
b}\Big(\sup\limits_{u-\frac{1}{n}\leq s\leq u}\mathbb{E}\big|x(s)-x(u-\frac{1}{n})\big|^2\Big)\Big\}=\delta,
\end{eqnarray*}
where $\delta$ is defined by (\ref{vsfde2-10}). Thus,
(\ref{vsfde-add-18}) shows that
$$
\mathbb{E}V\big(t,x(t)-x_n(t)\big) \leq
\mathbb{E}V\big(t_1,x(t_1)-x_n(t_1)\big),
$$
that is,
\begin{eqnarray*}
&&\mathbb{E}|x(t)-x_n(t)|^2\leq
e^{\int_{t_1}^t \sigma(u)du}\mathbb{E}|x(t_1)-x_n(t_1)|^2\\
&&+\int_{t_1}^t\varrho(u)e^{\int_{u}^t
\sigma(s)ds}du\Big(\varepsilon_n
+\sup\limits_{a-\tau_0\leq \theta\leq
t_1}\mathbb{E}|x(\theta)-x_n(\theta)|^2\Big)\\
&\leq
&\Big(e^{\int_{t_1}^t \sigma(u)du}
+\int_{t_1}^t\varrho(u)
e^{\int_u^t \sigma(s)ds}du\Big) \sup\limits_{a-\tau_0\leq
\theta\leq t_1}\mathbb{E}|x(\theta)-x_n(\theta)|^2+\Big(\int_{t_1}^t\varrho(u)
e^{\int_u^t \sigma(s)ds}du\Big)\varepsilon_n\\
&=& \Big(1
+\int_{t_1}^t\big(\sigma(u)+\varrho(u)\big)e^{\int_{u}^t
\sigma(s)ds}du\Big)\sup\limits_{a-\tau_0\leq
\theta\leq t_1}\mathbb{E}|x(\theta)-x_n(\theta)|^2+\Big(\int_{t_1}^t\varrho(u)
e^{\int_u^t \sigma(s)ds}du\Big)\varepsilon_n\\
&\leq & \big(1+\gamma_0\mu
e^{\alpha_0(b-a)}\big)\sup\limits_{a-\tau_0\leq \theta\leq
a+(i-1)\mu}\mathbb{E}|x(\theta)-x_n(\theta)|^2+\beta_0\mu
e^{\alpha_0(b-a)}\varepsilon_n
\end{eqnarray*}
for all $ t\in [a+(i-1)\mu, a+i\mu],\ \ i=1,2,\cdots,q$. Consequently,
\begin{eqnarray*}
&& \sup\limits_{a-\tau_0\leq
\theta\leq a+i\mu}\mathbb{E}|x(\theta)-x_n(\theta)|^2\\
&=&\max\Big\{\sup\limits_{a-\tau_0\leq \theta\leq
a+(i-1)\mu}\mathbb{E}|x(\theta)-x_n(\theta)|^2,
\sup\limits_{a+(i-1)\mu\leq
\theta\leq a+i\mu}\mathbb{E}|x(\theta)-x_n(\theta)|^2\Big\}\\
&\leq & \big(1+\gamma_0\mu
e^{\alpha_0(b-a)}\big)\sup\limits_{a-\tau_0\leq \theta\leq
a+(i-1)\mu}\mathbb{E}|x(\theta)-x_n(\theta)|^2+\beta_0\mu
e^{\alpha_0(b-a)}\varepsilon_n
\end{eqnarray*}
for $i=1,2,\cdots,q$. Therefore,
\begin{eqnarray}
&&\sup\limits_{a-\tau_0\leq \theta\leq
b}\mathbb{E}|x(\theta)-x_n(\theta)|^2=\sup\limits_{a-\tau_0\leq \theta\leq
a+q\mu}\mathbb{E}|x(\theta)-x_n(\theta)|^2\nonumber\\
&&\leq C_{\mu}^q \sup\limits_{a-\tau_0\leq \theta\leq
a}\mathbb{E}|x(\theta)-x_n(\theta)|^2+\frac{C_{\mu}^q-1}{\gamma_0\mu
e^{\alpha_0(b-a)}}\beta_0\mu e^{\alpha_0(b-a)}\varepsilon_n\label{vsfde2-13}\\
&&=\frac{\beta_0}{\gamma_0}(C_{\mu}^q-1)\varepsilon_n,\nonumber
\end{eqnarray}
where $C_{\mu}=1+\gamma_0\mu e^{\alpha_0(b-a)}.$
By Lemma \ref{lem-add-2}, we have
$
\varepsilon_n=\sup\limits_{a\leq t\leq
b}\sup\limits_{t-\frac{1}{n}\leq \theta\leq
t}\mathbb{E}\big|x(\theta)-x(t-\frac{1}{n})\big|^2\rightarrow 0,\
\ \mbox{as }n\rightarrow \infty.
$
Let $n\rightarrow \infty$ and take into account that
$$
\lim\limits_{q\rightarrow
\infty}C_{\mu}^q=\lim\limits_{q\rightarrow
\infty}\left(1+\frac{\gamma_0(b-a)e^{\alpha_0(b-a)}}{q}\right)^q=e^{\gamma_0(b-a)e^{\alpha_0(b-a)}}.
$$
Then (\ref{vsfde2-13}) leads to the relation (\ref{vsfde2-8}). The proof is complete.

\begin{thm}\label{thm1}
Assume that problem (\ref{vsfde1}) $\in
\mathcal{SD}(\alpha,\beta,\gamma_1,\gamma_2)$. Let
$c=\max\limits_{a\leq t\leq
b}$$\big(2\alpha(t)+2\beta(t)+\gamma_1^2(t)+2\gamma_1(t)\gamma_2(t)+\gamma_2^2(t)\big)$.
Then $\forall t\in [a,b]$,
\begin{eqnarray}
&\mathbb{E}|x(t)-y(t)|^2\leq e^{c(t-a)}\sup\limits_{a-\tau_0\leq \theta\leq
a}\mathbb{E}|\xi(\theta)-\eta(\theta)|^2, & \mbox{ if }c\geq 0, \label{vsfde2-15}\\
&\mathbb{E}|x(t)-y(t)|^2\leq \sup\limits_{a-\tau_0\leq \theta\leq
a}\mathbb{E}|\xi(\theta)-\eta(\theta)|^2,& \mbox{ if }c\leq 0,\label{vsfde-add-59}
\end{eqnarray}
where $y(t)$ is the solution of the perturbed problem (\ref{vsfde2-2}).
\end{thm}

The inequalities (\ref{vsfde2-15})
and (\ref{vsfde-add-59}) mean that problem (\ref{vsfde1}) is  stable in mean square and contractive in mean square, respectively.

\noindent{\bf Proof.}
We divide the proof into two cases: $\mu_1^{(0)}>0$ and $\mu_1^{(0)}=0$.

{\bf Case A: $\mu_1^{(0)}>0$}.
In this case, we can obtain the desired result in a similar manner
as in the proof of Theorem 2.1 in \cite{Lishoufu2005}. In fact,
replacing $\alpha(t), \beta(t), \|y(t)-z(t)\|$ in
\cite{Lishoufu2005} with
$2\alpha(t)+\beta(t)+\gamma_1(t)\gamma_2(t)+\gamma_1^2(t),
\beta(t)+\gamma_1(t)\gamma_2(t)+\gamma_2^2(t),
\mathbb{E}|x(t)-y(t)|^2$, respectively, using Lemma \ref{lem1} and
following the proof of Theorem 2.1 in \cite{Lishoufu2005}, we can
obtain either (\ref{vsfde2-15}) or (\ref{vsfde-add-59}) immediately.

{\bf Case B: $\mu_1^{(0)}=0$}.
Note that
\begin{eqnarray*}
\mathbb{E}|x(t)-y(t)|^2&\leq& \mathbb{E}\big(|x(t)-x_n(t)|+|x_n(t)-y_n(t)|+|y(t)-y_n(t)|\big)^2  \\
&\leq&3\left(\mathbb{E}|x(t)-x_n(t)|^2+\mathbb{E}|x_n(t)-y_n(t)|^2+\mathbb{E}|y(t)-y_n(t)|^2\right),
\end{eqnarray*}
where $y_n(t)$ is the solution of

\begin{numcases} {\label{vsfde-add-61}}
dy(t)=f(t,y(t),y^{(n,t)}(t-\tau(t)))dt + g(t,y(t),y^{(n,t)}(t-\tau(t)))dw(t), t \in[a,b], \quad \qquad \label{vsfde-add-61a}\\
y(t)=\eta(t), \, t\in [a-\tau_0,a], \ \ \eta\in
L^p_{\mathcal {F}_a} ([a-\tau_0, a];\mathbb{C}^d),\label{vsfde-add-61b}
\end{numcases}
which is the $\frac{1}{n}$-perturbed problem of (\ref{vsfde2-2}), and $y^{(n,t)}(t-\tau(t))$ is defined by
$$
y^{(n,t)}(t-\tau(t))=\left\{
\begin{array}{ll}
y(t-\tau(t)),& \tau(t)\geq \frac{1}{n},\\
y(t-\frac{1}{n}),&\tau(t)<\frac{1}{n}.
\end{array}\right.
$$
It is known that, problem (\ref{vsfde1}) $\in
\mathcal{SD}(\alpha,\beta,\gamma_1,\gamma_2)$ implies that problem
(\ref{vsfde10}) $\in \mathcal{SD}(\alpha$, $\beta$, $\gamma_1$,
$\gamma_2)$. It follows from (\ref{vsfde12}) that
$\tilde{\mu}_1^{(0)}\geq \frac{1}{n}>0$. Therefore, by case A, for
$\mathbb{E}|x_n(t)-y_n(t)|^2$, either (\ref{vsfde2-15}) holds if
$c>0$ or (\ref{vsfde-add-59}) holds if $c\leq 0$.
Letting $n\rightarrow +\infty$  and using Lemma \ref{lem2}, we can obtain the desired estimate of $\mathbb{E}|x(t)-y(t)|^2$ in this case.

\begin{cor}\label{cor-add-7}
Under the assumptions of Theorem \ref{thm1}.  Suppose $f(t,0,0)=0$
and $g(t,0,0)=0$, then $\forall t\in [a,b]$,
\begin{eqnarray}
&\mathbb{E}|x(t)|^2\leq e^{c(t-a)}\sup\limits_{a-\tau_0\leq \theta\leq
a}\mathbb{E}|\xi(\theta)|^2, & \mbox{ if }c>0, \nonumber\\
&\mathbb{E}|x(t)|^2\leq \sup\limits_{a-\tau_0\leq \theta\leq
a}\mathbb{E}|\xi(\theta)|^2,& \mbox{ if }c\leq 0.\nonumber
\end{eqnarray}
\end{cor}

\begin{lem}\label{lem-add-5}
Suppose problem (\ref{vsfde1}) $\in
\mathcal{SD}(\alpha,\beta,\gamma_1,\gamma_2)$, and that
\begin{equation}
\left.
\begin{array}{rl}
&2\alpha(t)+\beta(t)+\gamma_1(t)\gamma_2(t)+\gamma_1^2(t)\leq
\alpha_0<0,\\
&\frac{\beta(t)+\gamma_1(t)\gamma_2(t)+\gamma_2^2(t)}{|2\alpha(t)+\beta(t)+\gamma_1(t)\gamma_2(t)+\gamma_1^2(t)|}\leq
\nu<1,
\end{array}\right. \forall t\in [a,b],
\end{equation}
where $\alpha_0$ and $\nu$ are constants. Then for any given $c_1,
c_2, c_3: a\leq c_1<c_2<c_3\leq b$, we have
\begin{equation}\label{vsfde-add-81}
\begin{array}{rcl}
\quad \mathbb{E}|x(t)-y(t)|^2&\leq &
\Big(\nu+(1-\nu)e^{\alpha_0(c_2-c_1)}\Big)\sup\limits_{\mu_2^{(0)}(c_1,c_3)\leq
\theta\leq c_2}\mathbb{E}|x(\theta)-y(\theta)|^2,\\
&&\qquad\qquad \forall t\in
[c_2,c_3].
\end{array}
\end{equation}
\end{lem}

\noindent{\bf Proof.}
We divide the proof into two cases: $\mu_1^{(0)}>0$ and $\mu_1^{(0)}=0$.

{\bf Case A: $\mu_1^{(0)}>0$.}
In this case,  replacing $\alpha(t), \beta(t), \|y(t)-z(t)\|$ in
\cite{Lishoufu2005} with
$2\alpha(t)+\beta(t)+\gamma_1(t)\gamma_2(t)+\gamma_1^2(t)$,
$\beta(t)+\gamma_1(t)\gamma_2(t)+\gamma_2^2(t),
\mathbb{E}|x(t)-y(t)|^2$, respectively, using Lemma \ref{lem1} and
following the proof of Lemma 2.3 in \cite{Lishoufu2005},  we can
obtain (\ref{vsfde-add-81}) immediately.

{\bf Case B: $\mu_1^{(0)}=0$.}
In this case, we can obtain the estimate (\ref{vsfde-add-81}) in a
similar manner as in the proof of Case B of Theorem \ref{thm1}.

\subsection{Infinite interval}

Let us now proceed to discuss the equation (\ref{vsfde1}) which
satisfies conditions (\ref{vsfde-add-21})-(\ref{vsfde5}) but the
integration interval $[a,b]$ replaced by $[a,+\infty)$.
Accordingly,  interval $[a-\tau_0,b]$ is replaced by
$[a-\tau_0,+\infty)$, and the symbol $\mathcal{SD}(\alpha,\beta,\gamma_1,\gamma_2)$ is replaced by
$\mathcal{\overline{SD}}(\alpha,\beta,\gamma_1,\gamma_2)$.

\begin{thm}\label{thm6}
Assume that problem (\ref{vsfde1}) $\in
\mathcal{\overline{SD}}(\alpha,\beta,\gamma_1,\gamma_2)$, and
\begin{eqnarray*}
&&\lim\limits_{t\rightarrow +\infty}\big(t-\tau(t)\big)=+\infty,\ \
\sup\limits_{a\leq
t<+\infty}\big(2\alpha(t)+\beta(t)+\gamma_1(t)\gamma_2(t)+\gamma_1^2(t)\big)=\alpha_0<0,\\
&&\sup\limits_{a\leq
t<+\infty}\frac{\beta(t)+\gamma_1(t)\gamma_2(t)+\gamma_2^2(t)}{|2\alpha(t)+\beta(t)+\gamma_1(t)\gamma_2(t)
+\gamma_1^2(t)|}=\nu<1.
\end{eqnarray*}
Then,
 for any given constant $\mu>0$, there exists a strictly
increased sequence $\{t_k\}$ which diverges to $+\infty$ as
$k\rightarrow +\infty$, where $t_0=a$, such that
\begin{equation}\label{vsfde2-20}
\sup\limits_{t_k\leq t\leq t_{k+1}}\mathbb{E}|x(t)-y(t)|^2\leq
C_{\mu}^{k+1}\sup\limits_{a-\tau_0\leq t\leq
a}\mathbb{E}|\xi(t)-\eta(t)|^2, \ \ k=0,1,2,\cdots,
\end{equation}
where $C_{\mu}=\nu+(1-\nu)e^{\alpha_0\mu}\in (0,1)$. Hence,
\begin{equation}\label{vsfde2-22}
\lim\limits_{t\rightarrow +\infty}\mathbb{E}|x(t)-y(t)|^2=0.
\end{equation}
\end{thm}

\noindent{\bf Proof.}
It is obvious that (\ref{vsfde2-20}) implies
(\ref{vsfde2-22}). So, only the proof of (\ref{vsfde2-20}) is
required. First we construct a sequence $\{t_k\}$ by induction. Let
$t_0=a$. Suppose that $t_k$ is chosen appropriately, where $k\geq
0$. Because $\lim\limits_{t\rightarrow
+\infty}(t-\tau(t))=+\infty$, there exists a $M$ such that for all
$t\geq M$,  we have $t-\tau(t)\geq t_k$ and therefore
$\mu_2^{(0)}(M,+\infty)\geq t_k$. So we can choose $t_{k+1}=M+\mu$
and have the relation
\begin{equation}\label{vsfde-add-82}
t_k\leq \mu_2^{(0)}(t_{k+1}-\mu,+\infty)\leq t_{k+1}-\mu<t_{k+1}.
\end{equation}
Using (\ref{vsfde-add-82}) and Lemma \ref{lem-add-5}, we get
\begin{eqnarray*}
&&\sup\limits_{t_{k}\leq t\leq t_{k+1}}\mathbb{E}|x(t)-y(t)|^2\leq \big(\nu+(1-\nu)e^{\alpha_0\mu}\big)\sup\limits_{\mu_2^{(0)}(t_{k}-\mu,t_{k+1})\leq t\leq t_{k}}\mathbb{E}|x(t)-y(t)|^2\\
&\leq & C_{\mu}\sup\limits_{t_{k-1}\leq t\leq t_{k}}\mathbb{E}|x(t)-y(t)|^2
\leq  \ldots\leq C_{\mu}^{k+1}\sup\limits_{a-\tau_0\leq t\leq a}\mathbb{E}|\xi(t)-\eta(t)|^2.
\end{eqnarray*}
The proof is complete.

\begin{cor}\label{cor-add-6}
Under the same conditions as Theorem \ref{thm6}. Furthermore,
suppose that $f(t,0,0)=0, g(t,0,0)=0$, then
\begin{eqnarray}
&&\sup\limits_{t_k\leq t\leq t_{k+1}}\mathbb{E}|x(t)|^2\leq
C_{\mu}^{k+1}\sup\limits_{a-\tau_0\leq t\leq a}\mathbb{E}|\xi(t)|^2,
\ \ k=0,1,2,\cdots,\label{vsfde2-23}\\
&&\lim\limits_{t\rightarrow +\infty}\mathbb{E}|x(t)|^2=0.\label{vsfde2-24}
\end{eqnarray}
\end{cor}
\begin{rem}
Li \cite{Lishoufu2005} discussed the stability of nonlinear stiff
Volterra functional differential equations in Banach spaces.
Theorem \ref{thm1} and Theorem \ref{thm6} can be regarded as
generalizations of Theorem 2.1 and Theorem 2.2 of
\cite{Lishoufu2005} restricted in finite-dimensional Hilbert
spaces $\mathbb{C}^d$ and finitely many delays to the stochastic
version, respectively. It should be pointed out that the drift coefficient $f$
is required to be locally Lipschitz continuous in this paper,
whereas the condition is not required in
\cite{Lishoufu2005}. It is known that local Lipschitz continuity
is not a strong restriction.
\end{rem}

\subsection{Examples}

System (\ref{vsfde1}) includes the following three classes of
SDDEs as special cases

\begin{itemize}
\item SDDEs with constant delays:  $\tau(t)\equiv \tau$.
\end{itemize}

\begin{itemize}
\item Stochastic pantograph equations: $t-\tau(t)=qt$, where
$0<q<1$ is a constant.
\end{itemize}

\begin{itemize}
\item SDDEs with piecewise constant arguments: $t-\tau(t)=\lfloor
t-i\rfloor$, where $\lfloor t \rfloor$ denotes the largest integer
number less than or equal to $t$, $i$ is a nonnegative integer.
\end{itemize}
Therefore, Theorem \ref{thm1}, Theorem \ref{thm6}, Corollary
\ref{cor-add-7} and Corollary \ref{cor-add-6} are valid for the
three classes of SDDEs mentioned above.

\begin{exam}\label{example1}
Consider the linear SDDEs
\begin{equation}\label{vsfde-add-89}
\begin{array}{ll}
&dx(t)=\big(A_1(t)x(t)+A_2(t)x(t-\tau(t))+F(t)\big)dt\\
&\qquad\qquad+\big(B_1(t)x(t)+B_2(t)x(t-\tau(t))+G(t)\big)dw(t),
\end{array}
\end{equation}
where $A_1(t), A_2(t), B_1(t), B_2(t)\in \mathbb{C}^{d\times d},
F(t), G(t)\in \mathbb{C}^d$ are continuous with respect to $t$,
$w(t)$ is an 1-dimensional Wiener process. For the
problems (\ref{vsfde-add-89}),  it is easy to
verify the conditions (\ref{vsfde-add-21})-(\ref{vsfde5}) with
\begin{eqnarray*}
\alpha(t)=\lambda_{\max}^{\frac{A_1^*(t)+A_1(t)}{2}},\
\beta(t)=|A_2(t)|,\ \gamma_1(t)=|B_1(t)|,\ \gamma_2(t)=|B_2(t)|,
\end{eqnarray*}
where $\lambda_{\max}^{\frac{A_1^*(t)+A_1(t)}{2}}$ denotes the
largest eigenvalue of the Hermite matrix
$\frac{A_1^*(t)+A_1(t)}{2}$.
\end{exam}
Applying Theorem \ref{thm1} and Corollary \ref{cor-add-7} to
(\ref{vsfde-add-89}), we have the following corollary.
\begin{cor}\label{cor-add-8}
The solutions of (\ref{vsfde-add-89}) satisfy
\begin{eqnarray*}
&\mathbb{E}|x(t)-y(t)|^2\leq e^{c(t-a)}\sup\limits_{a-\tau_0\leq
\theta\leq
a}\mathbb{E}|\xi(\theta)-\eta(\theta)|^2,&t\in [a,b],  \mbox{ if }c\geq 0, \\
&\mathbb{E}|x(t)-y(t)|^2\leq \sup\limits_{a-\tau_0\leq \theta\leq
a}\mathbb{E}|\xi(\theta)-\eta(\theta)|^2,& t\in [a,b],  \mbox{ if
}c\leq 0,
\end{eqnarray*}
where $x(t), y(t)$ are the solutions of (\ref{vsfde-add-89}) corresponding to the
initial functions $\xi(t)$ and $\eta(t)$, respectively,
$$
c=\max\limits_{a\leq t\leq
b}\Big(2\lambda_{\max}^{\frac{A_1^*(t)+A_1(t)}{2}}+2|A_2(t)|+\big(|B_1(t)|+|B_2(t)|\big)^2\Big).
$$
Furthermore, if $F(t)=0$ and $G(t)=0$, then the solutions of
(\ref{vsfde-add-89}) satisfy
\begin{eqnarray*}
&\mathbb{E}|x(t)|^2\leq e^{c(t-a)}\sup\limits_{a-\tau_0\leq
\theta\leq
a}\mathbb{E}|\xi(\theta)|^2, & t\in [a,b],  \mbox{ if }c>0, \\
&\mathbb{E}|x(t)|^2\leq \sup\limits_{a-\tau_0\leq \theta\leq
a}\mathbb{E}|\xi(\theta)|^2,& t\in [a,b],  \mbox{ if }c\leq 0.
\end{eqnarray*}
\end{cor}
Applying Theorem \ref {thm6} and Corollary \ref {cor-add-6} to
(\ref{vsfde-add-89}) leads to the following
\begin{cor}\label{cor-add-9}
If $\lim\limits_{t\rightarrow +\infty}(t-\tau(t))=+\infty$,
\begin{eqnarray*}
&\sup\limits_{a\leq
t<+\infty}\Big(2\lambda_{\max}^{\frac{A_1^*(t)+A_1(t)}{2}}+|A_2(t)|+|B_1(t)||B_2(t)|+|B_1(t)|^2\Big)<0,\\
&\sup\limits_{a\leq
t<+\infty}\frac{|A_2(t)|+|B_1(t)||B_2(t)|+|B_2(t)|^2}{\big|2\lambda_{\max}^{\frac{A_1^*(t)+A_1(t)}{2}}
+|A_2(t)|+|B_1(t)||B_2(t)| +|B_1(t)|^2\big|}<1,
\end{eqnarray*}
then the solutions of (\ref{vsfde-add-89}) satisfy $\lim\limits_{t\rightarrow +\infty}\mathbb{E}|x(t)-y(t)|^2=0$.
Furthermore, if $F(t)=0$ and $G(t)=0$, then the solutions of
(\ref{vsfde-add-89}) satisfy $\lim\limits_{t\rightarrow +\infty}\mathbb{E}|x(t)|^2=0$.
\end{cor}
In particular, if $d=1$, $A_1(t), A_2(t), B_1(t), B_2(t)$ are
constants, that is, $A_1(t)=A_1, A_2(t)=A_2, B_1(t)=B_1,
B_2(t)=B_2$,  then the solutions of (\ref{vsfde-add-89}) satisfy $\lim\limits_{t\rightarrow +\infty}\mathbb{E}|x(t)|^2=0$
if
\begin{equation}\label{eqn-add-9}
\Re A_1+|A_2|+\frac{1}{2}\Big(|B_1|+|B_2|\Big)^2<0.
\end{equation}

\begin{rem}
For linear scalar SDDEs with constant delay and linear scalar
stochastic pantograph equations,  the condition (\ref{eqn-add-9}) is stated in \cite{Mao97,LCF04} and
\cite{Fan-Song-Liu-09}, respectively. Therefore, Corollary \ref{cor-add-6} includes as special cases the related results in \cite{Mao97,LCF04,Fan-Song-Liu-09}.
\end{rem}

\begin{exam}\label{example2} Consider the nonlinear equation
\begin{equation}\label{example-eqn-1}
\begin{array}{ll}
&dx(t)=\left( A_1(t)x+A_2(t)x^3+A_3(t)\sqrt{x^2(t-\tau(t))+1} +F(t)\right)dt \\
&\qquad\quad+\Big( B_1(t)\sin x(t)+B_2(t)\arctan x(t-\tau(t)) +G(t)\Big)dw(t),
\end{array}
\end{equation}
where $A_1(t), A_2(t), A_3(t), B_1(t), B_2(t), F(t), G(t)$ are continuous real-valued functions in $t$ and $A_2(t)<0$.
It is easy to
verify that (\ref{example-eqn-1}) satisfies the conditions (\ref{vsfde-add-21})-(\ref{vsfde5}) with
$$
\alpha(t)=A_1(t),\ \ \beta(t)=|A_3(t)|,\ \ \gamma_1(t)=|B_1(t)|,\ \ \gamma_2(t)=|B_2(t)|.
$$
\end{exam}
Applying Theorem \ref{thm1}, Corollary \ref{cor-add-7}, Theorem \ref {thm6} and Corollary \ref {cor-add-6} to (\ref{example-eqn-1}), we can derive the results of the solutions of (\ref{example-eqn-1}).
For the sake of brevity, we do not present them here.


\begin{rem}
The drift coefficient $f$ of (\ref{example-eqn-1}) satisfies local Lipschitz condition and one-sided Lipschitz condition with respect to $x$ but global Lipschitz condition. The stability analysis in this work is based on the local Lipschitz condition and the one-sided Lipschitz condition, rather than a more restrictive global Lipschitz condition.
\end{rem}

\section{Stability of backward Euler method}

In this section, we investigate whether numerical methods can
reproduce the contractivity in mean square. For the deterministic
differential equations, it is known that the contractivity of numerical methods is too strong \cite{BZ92,HLFC02}. The existing
theories \cite{BZ92,Torelli89} show that only the backward Euler method
and the two-stage Lobatto IIIC method can preserve the
contractivity of nonlinear delay differential equations.
Therefore, in the stochastic setting, we only focus on the
backward Euler method instead of other methods.

For simplicity, from now on, we assume that
$$
\alpha(t)\equiv \alpha,\ \beta(t)\equiv \beta,\ \gamma_1(t)\equiv\gamma_1,\ \gamma_2(t)\equiv\gamma_2, \ \ t\in [a,b].
$$

\noindent On a finite time interval $[a,b]$, a uniformly partition is defined by
$$
t_i=a+ih, \ \ i=0,1,\ldots, h=\frac{b-a}{N}.
$$
The backward Euler method applied to (\ref{vsfde1}) yields
\begin{numcases}{\label{vsfde-add-36}}
\left.\begin{array}{rcl}
X_{n+1}&=&X_n+hf(t_{n+1},X_{n+1},X^h(t_{n+1}-\tau(t_{n+1})))\\
&&\ +g(t_n,X_n,X^h(t_n-\tau(t_n)))
\Delta w_n,\quad n=0,1,\ldots, N-1,\\
\end{array}\right.
\label{vsfde-add-36-1}\\
X^h(t)=\pi^h(t,\xi, X_1,X_2,\ldots,X_{n}),\ \ a-\tau_0\leq t\leq
t_{n},\label{vsfde-add-36-2}
\end{numcases}
where $\pi^h$ is an appropriate interpolation operator which approximates to the exact solution $x(t)$ on the interval $[a-\tau_0,b]$, $X_n$ is an approximation to the exact solution $x(t_n)$,  $\Delta w_n=w(t_{n+1})-w(t_n)$.
It is well known that the backward Euler method is convergent with strong order only 1/2
for stochastic differential equations. So, interpolation operator $\pi^h$ could be chosen as the follows
\begin{equation}\label{vsfde-add-37}
 X^h(t)=\left\{
\begin{array}{l}
\frac{1}{h} [(t_{i+1}-t)X_i+(t-t_i)X_{i+1}], t_i\leq t\leq t_{i+1},  i=0,1,2,\ldots,N-1,\\
\xi(t),\ \ a-\tau_0\leq t\leq a.
\end{array}\right.
\end{equation}
Applying the backward Euler method to the perturbed problem (\ref{vsfde2-2}) we can obtain the corresponding scheme
\begin{equation}{\label{vsfde-add-38}}
\left\{\begin{array}{l}
Y_{n+1}=Y_n+hf(t_{n+1},Y_{n+1},Y^h(t_{n+1}-\tau(t_{n+1})))\\
\qquad\qquad +g(t_n,Y_n,Y^h(t_n-\tau(t_n)))\Delta w_n,\qquad \ n=0,1,\ldots, N-1,\\
Y^h(t)=\pi^h(t,\eta, Y_1,Y_2,\ldots,Y_{n}),\ \ a-\tau_0\leq t\leq
t_{n}.
\end{array}\right.
\end{equation}
For simplicity, for any given nonnegative integer $n$, we write
\begin{equation}\label{abbreviation1}
\begin{array}{lll}
&&P_n=X_n-Y_n,\ \  Q_n=\max\{\max\limits_{1\leq i\leq n}\mathbb{E}|P_i|^2, \sup\limits_{a-\tau_0\leq t\leq a}\mathbb{E}|\xi(t)-\eta(t)|^2\},\ \ n\geq 1,\\
&&Q_0=\sup\limits_{a-\tau_0\leq t\leq a}\mathbb{E}|\xi(t)-\eta(t)|^2.
\end{array}
\end{equation}
Moreover, for convenience, we introduce notations to denote the values of drift and diffusion coefficients at specific points.
\begin{equation}\label{abbreviation2}
\begin{array}{ll}
&f^{xx}(n+1)=f(t_{n+1},X_{n+1},X^h(t_{n+1}-\tau(t_{n+1}))),\\
&f^{yy}(n+1)=f(t_{n+1},Y_{n+1},Y^h(t_{n+1}-\tau(t_{n+1}))),\\
&f^{yx}(n+1)=f(t_{n+1},Y_{n+1},X^h(t_{n+1}-\tau(t_{n+1}))),\\
&g^{xx}(n)=g(t_n,X_n,X^h(t_n-\tau(t_n))),\ g^{yy}(n)=g(t_n,Y_n,Y^h(t_n-\tau(t_n))).
\end{array}
\end{equation}

\begin{lem}\label{lem-add-6}
Under the conditions  (\ref{vsfde3}) and (\ref{vsfde4}), if $(\alpha+\beta)h<1$, the implicit equation (\ref{vsfde-add-36-1}) admits a  unique solution.
\end{lem}

\noindent{\bf Proof.}
Let $\tilde{f}(z)=f(\cdot,z, z^h(\cdot))$, then
implicit equation (\ref{vsfde-add-36-1}) can be rewritten as
\begin{equation}\label{vsfde-add-85}
z=h\tilde{f}(z)+b=hf(\cdot,z,z^h(\cdot))+b,
\end{equation}
where $z$ is unknown whereas $b$ and $h$ are known.
Inserting the interpolation operator (\ref{vsfde-add-37}) into
(\ref{vsfde-add-85}), we have
\begin{equation}\label{vsfde-add-86}
z=h\tilde{f}(z)+b=hf(\cdot,z,lz+b_0)+b,
\end{equation}
where $0\leq l\leq 1$, $l$ and $b_0$ are also known.
It follows from (\ref{vsfde3}), (\ref{vsfde4}) and
(\ref{vsfde-add-86}) that
\begin{eqnarray*}
&&\Re \langle z_1-z_2, \tilde{f}(z_1)-\tilde{f}(z_2) \rangle=\Re
\langle z_1-z_2, f(\cdot,z_1, lz_1+b_0)
-f(\cdot,z_2, lz_2+b_0)\rangle\\
&=&\Re \langle z_1-z_2, f(\cdot,z_1, lz_1+b_0)-f(\cdot,z_2,
lz_1+b_0)\rangle\\
&&+\Re \langle z_1-z_2, f(\cdot,z_2, lz_1+b_0)-f(\cdot,z_2, lz_2+b_0)\rangle
\leq  \alpha |z_1-z_2|^2+\beta |z_1-z_2|^2.
\end{eqnarray*}
The assertion follows immediately from Theorem 5.6.1 in \cite{HW96}.

\begin{thm}\label{thm-add-1}
Assume that problem (\ref{vsfde1}) $\in
\mathcal{SD}(\alpha,\beta,\gamma_1,\gamma_2)$. Let $\{X_n\}$ and
$\{Y_n\}$ be two sequences of numerical solutions obtained by the
backward Euler schemes (\ref{vsfde-add-36}) and
(\ref{vsfde-add-38}), respectively.  Write
$c=2\alpha+2\beta+\gamma_1^2+2\gamma_1\gamma_2+\gamma_2^2$.

(i) If $c>0$, for any given $c_0\in (0,1)$, then we have for $hc\leq c_0$
\begin{equation}\label{vsfde-add-41}
\mathbb{E}|X_n-Y_n|^2\leq e^{\tilde{c}(t_n-a)}\sup\limits_{a-\tau_0\leq t\leq a}\mathbb{E}|\xi(t)-\eta(t)|^2, \ \  \ n=1,2,\cdots,N,
\end{equation}
where $\tilde{c}=\frac{c_1}{h}$,
\begin{eqnarray*}
c_1&=&\max\left\{\frac{1+h\gamma_1^2+2h\gamma_1\gamma_2+h\gamma_2^2}{1-2h\alpha-2h\beta},
\frac{1+h\beta+h\gamma_1^2+2h\gamma_1\gamma_2+h\gamma_2^2}{1-2h\alpha-h\beta}\right\}\\
&=&\frac{1+h\gamma_1^2+2h\gamma_1\gamma_2+h\gamma_2^2}{1-2h\alpha-2h\beta}>1.
\end{eqnarray*}

(ii) If $c\leq 0$, then we have for any $h>0$
\begin{equation}\label{vsfde-add-42}
\mathbb{E}|X_n-Y_n|^2\leq \sup\limits_{a-\tau_0\leq t\leq a}\mathbb{E}|\xi(t)-\eta(t)|^2, \ \  \ n=1,2,\cdots,N.
\end{equation}
\end{thm}

Note that (\ref{vsfde-add-41}) and (\ref{vsfde-add-42}) can be regarded as numerical analogs of (\ref{vsfde2-15}) and (\ref{vsfde-add-59}), respectively.

\noindent{\bf Proof.} (i)  By (\ref{vsfde-add-36}) and (\ref{vsfde-add-38}), we have
\begin{eqnarray*}
P_{n+1}-h\big(f^{xx}(n+1)-f^{yy}(n+1)\big)=P_n+\big(g^{xx}(n)-g^{yy}(n)\big)\Delta
w_n,
\end{eqnarray*}
which yields
\begin{eqnarray*}
&&|P_{n+1}|^2-2h\Re \big\langle P_{n+1},
f^{xx}(n+1)-f^{yy}(n+1)
\big\rangle+h^2\big|f^{xx}(n+1)-f^{yy}(n+1)\big|^2\\
&=&|P_n|^2+2\Re \big\langle P_n, \big(g^{xx}(n)-g^{yy}(n)\big)\Delta w_n \big\rangle+\big|\big(g^{xx}(n)-g^{yy}(n)\big)\Delta
w_n\big|^2.
\end{eqnarray*}
Taking expectation and using (\ref{vsfde3})-(\ref{vsfde5}) and
(\ref{vsfde-add-37}), we get
\begin{equation}
\begin{array}{rcl}
\mathbb{E}|P_{n+1}|^2 & \leq & \mathbb{E}|P_n|^2+2h\mathbb{E}\Re
\big\langle P_{n+1}, f^{xx}(n+1)-f^{yy}(n+1) \big\rangle+h\mathbb{E}\big|g^{xx}(n)-g^{yy}(n)\big|^2\\
&\leq & \mathbb{E}|P_n|^2+2h\mathbb{E}\Re \big\langle P_{n+1},
f^{xx}(n+1)-f^{yx}(n+1) \big\rangle\\
&&+2h\mathbb{E}\Re \big\langle P_{n+1},
f^{yx}(n+1)-f^{yy}(n+1)
\big\rangle+h\mathbb{E}\big|g^{xx}(n)-g^{yy}(n)\big|^2\\
&\leq & \mathbb{E}|P_n|^2+2h\alpha \mathbb{E}|P_{n+1}|^2\\
&&+2h\beta
\mathbb{E}\big(|P_{n+1}||X^h(t_{n+1}-\tau(t_{n+1}))-Y^h(t_{n+1}-\tau(t_{n+1}))|\big)\\
&&+h\mathbb{E}\big(\gamma_1|P_n|+\gamma_2|X^h(t_n-\tau(t_n))-Y^h(t_n-\tau(t_n))|\big)^2\\
&\leq & \mathbb{E}|P_n|^2+2h\alpha \mathbb{E}|P_{n+1}|^2\\
&&+h\beta \mathbb{E}|P_{n+1}|^2+h\beta
\mathbb{E}|X^h(t_{n+1}-\tau(t_{n+1}))-Y^h(t_{n+1}-\tau(t_{n+1}))|^2\\
&&+h(\gamma_1^2+\gamma_1\gamma_2) \mathbb{E}|P_n|^2+h(\gamma_2^2+\gamma_1\gamma_2)
\mathbb{E}|X^h(t_n-\tau(t_n))-Y^h(t_n-\tau(t_n))|^2\\
&\leq & \mathbb{E}|P_n|^2+2h\alpha \mathbb{E}|P_{n+1}|^2+h\beta \mathbb{E}|P_{n+1}|^2\\
&&+h\beta \max\big\{\max\limits_{1\leq i\leq
n+1}\mathbb{E}|P_i|^2,
\sup\limits_{a-\tau_0\leq t\leq a}\mathbb{E}|\xi(t)-\eta(t)|^2\big\}\\
&&+h(\gamma_1^2+\gamma_1\gamma_2)
\mathbb{E}|P_n|^2+h(\gamma_2^2+\gamma_1\gamma_2)  Q_n,
\end{array}
\label{vsfde-add-65}
\end{equation}
where we used the piecewise linear interpolation
(\ref{vsfde-add-37}) and the following inequality
\begin{equation}\label{vsfde-add-88}
\mathbb{E}|(1-\delta)P_i+\delta P_{i+1}|^2\leq
\max\{\mathbb{E}|P_i|^2, \mathbb{E}|P_{i+1}|^2\}, \ \ 0\leq
\delta\leq 1.\\
\end{equation}
It is clear from (\ref{vsfde-add-65}) that
\begin{eqnarray}
&&(1-2h\alpha -h\beta)\mathbb{E}|P_{n+1}|^2\leq (1+h\gamma_1^2+h\gamma_1\gamma_2)\mathbb{E}|P_n|^2\nonumber
\\[-1.5ex]\label{vsfde-add-40}\\[-1.5ex]
&&\quad+h\beta\max\big\{\max\limits_{1\leq i\leq n+1}\mathbb{E}|P_i|^2,
\sup\limits_{a-\tau_0\leq t\leq
a}\mathbb{E}|\xi(t)-\eta(t)|^2\big\}+h(\gamma_2^2+\gamma_1\gamma_2)
Q_n,\nonumber
\end{eqnarray}
We now consider two cases:
\begin{eqnarray*}
\mbox{{\bf(a)} } \quad \max\{\max\limits_{1\leq i\leq n+1}\mathbb{E}|P_i|^2, \sup\limits_{a-\tau_0\leq t\leq a}\mathbb{E}|\xi(t)-\eta(t)|^2\}=\mathbb{E}|P_{n+1}|^2,\\
\mbox{{\bf(b)} } \quad \max\{\max\limits_{1\leq i\leq n+1}\mathbb{E}|P_i|^2, \sup\limits_{a-\tau_0\leq t\leq a}\mathbb{E}|\xi(t)-\eta(t)|^2\}\ne \mathbb{E}|P_{n+1}|^2.\\
\end{eqnarray*}
In the case of {\bf (a)}, it follows from (\ref{vsfde-add-40}) that
\begin{equation}
(1-2h\alpha-2h\beta)\mathbb{E}|P_{n+1}|^2\leq
(1+h\gamma_1^2+2h\gamma_1\gamma_2+h\gamma_2^2)Q_n,
\end{equation}
which yields
\begin{equation}\label{vsfde-add-43}
\mathbb{E}|P_{n+1}|^2\leq
\frac{1+h\gamma_1^2+2h\gamma_1\gamma_2+h\gamma_2^2}{1-2h\alpha-2h\beta}Q_n\leq
c_1Q_n.
\end{equation}
In the case of {\bf (b)}, (\ref{vsfde-add-40}) implies that
\begin{eqnarray*}
&&(1-2h\alpha-h\beta)\mathbb{E}|P_{n+1}|^2\leq (1+h\gamma_1^2+h\gamma_1\gamma_2)\mathbb{E}|P_n|^2
+h\beta Q_n+h(\gamma_2^2+\gamma_1\gamma_2)Q_n\\
&\leq & (1+h\beta+h\gamma_1^2+2h\gamma_1\gamma_2+h\gamma_2^2)Q_n,
\end{eqnarray*}
which yields
\begin{equation}\label{vsfde-add-44}
\mathbb{E}|P_{n+1}|^2\leq
\frac{1+h\beta+h\gamma_1^2+2h\gamma_1\gamma_2+h\gamma_2^2}{1-2h\alpha-h\beta}Q_n\leq
c_1Q_n.
\end{equation}
To summarize, both in the cases we have shown that $\mathbb{E}|P_{n+1}|^2\leq c_1Q_n$,
which yields
\begin{equation}\label{vsfde-add-45}
Q_n\leq Q_{n-1}+\mathbb{E}|P_n|^2\leq (1+c_1)Q_{n-1}.
\end{equation}
By induction, we further obtain
\begin{eqnarray*}
\mathbb{E}|X_n-Y_n|^2&=&\mathbb{E}|P_n|^2\leq Q_n\leq
(1+c_1)Q_{n-1}\leq \cdots \leq (1+c_1)^nQ_0\\
&\leq & e^{c_1n}Q_0=
e^{\tilde{c}(t_n-a)}\sup\limits_{a-\tau_0\leq t\leq
a}\mathbb{E}|\xi(t)-\eta(t)|^2.
\end{eqnarray*}

(ii) When $c\leq 0$, noting that (\ref{vsfde-add-43}), (\ref{vsfde-add-44}) and
$$
\frac{1+h\gamma_1^2+2h\gamma_1\gamma_2+h\gamma_2^2}{1-2h\alpha-2h\beta}\leq
1,\ \
\frac{1+h\beta+h\gamma_1^2+2h\gamma_1\gamma_2+h\gamma_2^2}{1-2h\alpha-h\beta}\leq
1,
$$
we have for any $h>0$
\begin{equation}
\mathbb{E}|X_n-Y_n|^2\leq Q_{n-1}\leq Q_{n-2}\leq \cdots\leq Q_0=\sup\limits_{a-\tau_0\leq t\leq a}\mathbb{E}|\xi(t)-\eta(t)|^2.
\end{equation}
Therefore we have completed the proof of the theorem.

\begin{cor}\label{cor-add-1}
Assume that problem (\ref{vsfde1}) $\in
\mathcal{SD}(\alpha,\beta,\gamma_1,\gamma_2)$. Let $\{X_n\}$ be a
sequence of numerical solutions obtained by the backward Euler
method (\ref{vsfde-add-36}).  Furthermore, if $f(t,0,0)=0,
g(t,0,0)=0$, and

(i) if $c>0$, for any given $c_0\in (0,1)$, then we have for $hc\leq c_0$
\begin{equation}\label{vsfde-add-92}
\mathbb{E}|X_n|^2\leq e^{\tilde{c}(t_n-a)}\sup\limits_{a-\tau_0\leq t\leq a}\mathbb{E}|\xi(t)|^2, \ \  n=1,2,\cdots,N;
\end{equation}

(ii) if $c\leq 0$, then we have for any $h>0$
\begin{equation}\label{vsfde-add-93}
\mathbb{E}|X_n|^2\leq \sup\limits_{a-\tau_0\leq t\leq a}\mathbb{E}|\xi(t)|^2, \ \  n=1,2,\cdots,N.
\end{equation}
\end{cor}

\begin{thm}\label{thm-add-2}
Assume that problem (\ref{vsfde1}) $\in
\mathcal{\overline{SD}}(\alpha,\beta,\gamma_1,\gamma_2)$, and
\begin{equation}\label{vsfde-add-94}
\lim\limits_{t\rightarrow +\infty}(t-\tau(t))=+\infty,\ \
c=2\alpha+2\beta+\gamma_1^2+2\gamma_1\gamma_2+\gamma_2^2<0.
\end{equation}
Let $\{X_n\}$ and $\{Y_n\}$ be two sequences of numerical solutions obtained by the backward Euler schemes (\ref{vsfde-add-36}) and (\ref{vsfde-add-38}).
Then,

(i)  there exists a strictly increased positive integer sequence $\{n_k\}$ which diverges to $+\infty$ as $k\rightarrow +\infty$, where
$n_0=0$, such that for any given $h>0$,
\begin{equation}\label{vsfde-add-46}
\max\limits_{n_k<i\leq n_{k+1}}\mathbb{E}|X_i-Y_i|^2\leq c_2^{k+1}\sup\limits_{a-\tau_0\leq t\leq a}\mathbb{E}|\xi(t)-\eta(t)|^2 ,\ \ k=0,1,2,\ldots,
\end{equation}
where
\begin{eqnarray*}
c_2&=&\max\left\{\frac{1+h\gamma_1^2+2h\gamma_1\gamma_2+h\gamma_2^2}{1-2h\alpha-2h\beta},
\frac{1+h\beta+h\gamma_1^2+2h\gamma_1\gamma_2+h\gamma_2^2}{1-2h\alpha-h\beta}\right\}\\
&=&\frac{1+h\beta+h\gamma_1^2+2h\gamma_1\gamma_2+h\gamma_2^2}{1-2h\alpha-h\beta}<1;
\end{eqnarray*}

(ii) for any given $h>0$,
\begin{equation}\label{vsfde-add-47}
\lim\limits_{n\rightarrow +\infty}\mathbb{E}|X_n-Y_n|^2=0.
\end{equation}
\end{thm}
Note that (\ref{vsfde-add-46})
and (\ref{vsfde-add-47}) can be regarded
as numerical analogs of (\ref{vsfde2-20}) and (\ref{vsfde2-22}), respectively.

\noindent{\bf Proof.}
It is obvious that (\ref{vsfde-add-46}) implies
(\ref{vsfde-add-47}), and we only need to prove
(\ref{vsfde-add-46}). By (\ref{vsfde-add-94}), we have
$2\alpha+\beta+\gamma_1\gamma_2+\gamma_1^2<0,
(\beta+\gamma_1\gamma_2+\gamma_2^2)/|2\alpha+\beta+\gamma_1\gamma_2+\gamma_1^2|<1$
and $c_2<1$.

First, as done in \cite{Li2007, WZ2010}, we can construct a strictly increased sequence of integers $\{n_k\}$ which diverges to $+\infty$ as $k \rightarrow +\infty$, such that
$$
t-\tau(t)>t_{n_k+1},\ \ \forall t\geq t_{n_{k+1}},
$$
where $n_0=0$. In fact, suppose that $n_k (k\geq 0)$ has been
chosen appropriately. Then there exists a constant $M>t_{n_k}$
such that for all $t\geq M$ we have $t-\tau(t)>t_{n_k}+h$ since
$\lim_{t\rightarrow +\infty}(t-\tau(t))=+\infty$. If $M$ is a
node, we let $t_{n_{k+1}}=M$, otherwise there exists natural
number $j$ such that $t_j<M<t_{j+1}$, then we let $n_{k+1}=j+1$
and $t_{n_{k+1}}=t_{j+1}$. Thus we obtain the required sequence
$\{n_k\}$ which satisfies
$$
t_0<t_1<\cdots<t_{n_1}< t_{n_1+1}<\cdots<t_{n_2}\cdots<t_{n_k}<\cdots.
$$

For $n_k<n+1\leq n_{k+1}$, by the second inequality of (\ref{vsfde-add-65}) and conditions (\ref{vsfde3})-(\ref{vsfde5}),
we have
\begin{eqnarray*}
\mathbb{E}|P_{n+1}|^2&\leq & \mathbb{E}|P_n|^2+2h\mathbb{E}\Re
\big\langle P_{n+1}, f^{xx}(n+1)
-f^{yx}(n+1) \big\rangle\nonumber\\
&&+2h\mathbb{E}\Re \big\langle P_{n+1},
f^{yx}(n+1)-f^{yy}(n+1)
\big\rangle+h\mathbb{E}\big|g^{xx}(n)-g^{yy}(n)\big|^2\nonumber\\
&\leq & \mathbb{E}|P_n|^2+2h \alpha
\mathbb{E}|P_{n+1}|^2+2h\beta\mathbb{E}\big(|P_{n+1}||X^h(t_{n+1}-\tau(t_{n+1}))-Y^h(t_{n+1}-\tau(t_{n+1}))|\big)\\
&&+h\mathbb{E}\Big(\gamma_1|P_n|+\gamma_2|X^h(t_n-\tau(t_n))-Y^h(t_n-\tau(t_n))| \Big)^2\\
&\leq & \mathbb{E}|P_n|^2+2h \alpha \mathbb{E}|P_{n+1}|^2+h\beta
\mathbb{E}|P_{n+1}|^2 +h\beta\mathbb{E}|X^h(t_{n+1}-\tau(t_{n+1}))-Y^h(t_{n+1}-\tau(t_{n+1}))|^2\\
&&+h(\gamma_1^2+\gamma_1\gamma_2)\mathbb{E}|P_n|^2+h(\gamma_2^2+\gamma_1\gamma_2)
\mathbb{E}|X^h(t_n-\tau(t_n))-Y^h(t_n-\tau(t_n))|^2,\\
\end{eqnarray*}
which yields
\begin{eqnarray*}
&&(1-2h\alpha-h\beta)\mathbb{E}|P_{n+1}|^2\\
&\leq &(1+h\gamma_1^2+h\gamma_1\gamma_2)\mathbb{E}|P_n|^2
+h\beta\mathbb{E}|X^h(t_{n+1}-\tau(t_{n+1}))-Y^h(t_{n+1}-\tau(t_{n+1}))|^2\\
&&+h(\gamma_2^2+\gamma_1\gamma_2)\mathbb{E}|X^h(t_n-\tau(t_n))-Y^h(t_n-\tau(t_n))|^2\\
&\leq&
(1+h\gamma_1^2+h\gamma_1\gamma_2)\mathbb{E}|P_n|^2+h\beta\max\limits_{n_{k-1}<
i\leq n+1}\mathbb{E}|P_i|^2
+h(\gamma_2^2+\gamma_1\gamma_2)\max\limits_{n_{k-1}< i\leq
n}\mathbb{E}|P_i|^2,
\end{eqnarray*}
where we used the piecewise linear interpolation operator
(\ref{vsfde-add-37}) and the inequality (\ref{vsfde-add-88}).
We now consider the following two  cases.

\noindent If $\max\limits_{n_{k-1}< i\leq n+1}\mathbb{E}|P_i|^2=\mathbb{E}|P_{n+1}|^2$,  we have
\begin{eqnarray*}
\mathbb{E}|P_{n+1}|^2\leq
\frac{1+h\gamma_1^2+2h\gamma_1\gamma_2+h\gamma_2^2}{1-2h\alpha-2h\beta}\max\limits_{n_{k-1}<
i\leq n}\mathbb{E}|P_i|^2\leq c_2\max\limits_{n_{k-1}< i\leq
n}\mathbb{E}|P_i|^2.
\end{eqnarray*}

\noindent If $\max\limits_{n_{k-1}< i\leq n+1}\mathbb{E}|P_i|^2\ne \mathbb{E}|P_{n+1}|^2$,  we have
\begin{eqnarray*}
\mathbb{E}|P_{n+1}|^2\leq
\frac{1+h\beta+h\gamma_1^2+2h\gamma_1\gamma_2+h\gamma_2^2}{1-2h\alpha-h\beta}\max\limits_{n_{k-1}<i\leq
n}\mathbb{E}|P_i|^2\leq c_2\max\limits_{n_{k-1}< i\leq
n}\mathbb{E}|P_i|^2.
\end{eqnarray*}
In both cases, we have
\begin{equation}\label{vsfde-add-52}
\mathbb{E}|P_{n+1}|^2\leq c_2\max\limits_{n_{k-1}< i\leq n}\mathbb{E}|P_i|^2,\ \ n_k<n+1\leq n_{k+1}.
\end{equation}
(\ref{vsfde-add-52}) with $n=n_k$ reduces to
$
\mathbb{E}|P_{n_k+1}|^2\leq c_2\max\limits_{n_{k-1}< i\leq n_k}\mathbb{E}|P_i|^2.
$
By induction, we have
\begin{eqnarray*}
\max\limits_{n_k<i\leq n_{k+1}}\mathbb{E}|X_i-Y_i|^2&=&\max\limits_{n_k<i\leq n_{k+1}}\mathbb{E}|P_i|^2\leq c_2\max\limits_{n_{k-1}< i\leq n_k}\mathbb{E}|P_i|^2
\\
&\leq & \cdots\leq c_2^{k+1}\max\limits_{a-\tau_0\leq t\leq a}\mathbb{E}|\xi(t)-\eta(t)|^2.
\end{eqnarray*}
The proof is complete.

\begin{cor}\label{cor-add-2}
Under the same assumptions of Theorem \ref{thm-add-2}. Let $\{X_n\}$ be a sequence of numerical solution obtained by the backward Euler method (\ref{vsfde-add-36}). Furthermore, if $f(t,0,0)=0$ and $g(t,0,0)=0$, then,

(i)  there exists a strictly increased positive integer sequence $\{n_k\}$ which diverges to $+\infty$ as $k\rightarrow +\infty$, where
$n_0=0$, such that for any given $h>0$,
\begin{eqnarray*}
\max\limits_{n_k<i\leq n_{k+1}}\mathbb{E}|X_i|^2\leq c_2^{k+1}\sup\limits_{a-\tau_0\leq t\leq a}\mathbb{E}|\xi(t)|^2 ,\ \ k=0,1,2,\ldots.
\end{eqnarray*}

(ii)  for any given $h>0,\
\lim\limits_{n\rightarrow +\infty}\mathbb{E}|X_n|^2=0.$
\end{cor}

\section{SDDEs with several delays}

Consider the following SDDEs with several delays
\begin{equation} {\label{eqn-add-7}}
\left\{\begin{array}{ll}
dx(t)=&f(t,x(t),x(t-\tau_1(t)),\cdots,x(t-\tau_r(t)))dt \\
&+ g(t,x(t),x(t-\tau_1(t)),\cdots,x(t-\tau_r(t)))dw(t),
 \quad  t\geq a,\\
x(t)=&\xi(t), \ \  t\in [a-\tau_0,a],
\end{array}\right.
\end{equation}
where $\tau_i(t)\geq 0, i=1,2,\cdots,r$ and $\max\limits_{1\leq
i\leq r}\inf\limits_{t\geq a}(t-\tau_i(t))\geq a-\tau_0$. All
results given in this paper can be extended easily to the case of
several delays. For the sake of brevity, we do not present the
corresponding results for (\ref{eqn-add-7}).

\section{Conclusions and future work}

In this paper, we investigate the stability of analytical and numerical solutions of
nonlinear SDDEs. We derive sufficient conditions for the stability, contractivity and asymptotic contractivity in mean square of
the solutions for nonlinear SDDEs. The results provide a unified
theoretical treatment for SDEs,
SDDEs with constant delay and variable
delay (including bounded and unbounded variable delays). Then, it is proved
that the backward Euler method can preserve the properties of the
underlying system.
The main results of analytic solution in this
paper can be regarded as a generalization of those in
\cite{Lishoufu2005} restricted in finite-dimensional Hilbert
spaces and finitely many delays to the stochastic version. We have
encountered problems when we tried to obtain a unified framework
for general SFDEs. It is
worth noting that whether the results in \cite{Lishoufu2005} can
be extended to general SFDEs or not. One area for the future work is to give a positive or
negative answer for the question.
Neutral stochastic delay differential equation (NSDDE) is more
general than stochastic delay differential equation. It is
interesting to investigate whether the theory of this paper can be
extended to NSDDEs and corresponding numerical methods. It will
also be our future work.


\begin{thebibliography}{10}

\bibitem{Baker-Buckwar05}  C.T.H.Baker and E.Buckwar, Exponential stability in $p$-th mean of
solutions, and of convergent Euler-type solutions, of stochastic
delay differential equations, J Comput. Appl. Math.
184 (2005) 404-427.

\bibitem{BZ92}A.Bellen, M.Zennaro, Strong contractivity properties of numerical methods for
ordinary and delay differential equations, Appl.Numer.Math.
9 (1992) 321-346.


\bibitem{Fan-Song-Liu-09}Z.Fan, M.Song, M.Liu,   The $\alpha$th moment
stability for the stochastic pantograph equation, J.Comput. Appl.
Math. 233 (2009) 109-120.

\bibitem{HW96}E.Hairer, G.Wanner, Solving Ordinary Differential Equations II: Stiff and Differential-
Algebraic Problems, Second ed., Springer-Verlag, Berlin, 1996.

\bibitem{HMS2002}D.J.Higham, X.Mao, A.M.Stuart, Strong convergence of Euler-type methods for
nonlinear stochastic differential equations, SIAM J. Numer. Anal.
40 (2002) 1041-1063.

\bibitem{Higham-Kloeden-2005}D.J.Higham, P.E.Kloeden, Numerical methods for nonlinear stochastic
differential equations with jumps, Numer.Math.  101 (2005) 101-119.

\bibitem{HLFC02}C.Huang, S.Li, H.Fu, G.Chen,  Nonlinear stability of general
linear methods for delay differential equations, BIT.
42 (2002) 380-392.

\bibitem{Lishoufu2005}S.Li, Stability analysis of solutions to nonlinear stiff
Volterra functional differential equations in Banach spaces,
Sci.China (Ser. A). 48 (2005) 372-387.

\bibitem{Li2007}S.Li, Contractivity and asymptotic stability properties of Runge-Kutta methods for Volterra functional
differential equations. In: International Workshop on Numerical
Analysis and Computational Methods for Functional Differential and
Integral Equations, Hong Kong Baptist University, 2007.


\bibitem{LCF04}M.Liu, W.Cao, Z.Fan,  Convergence and stability of the semi-implicit euler method for a linear stochastic
differential delay equation, J. Comput. Appl. Math.
170 (2004) 255-268.

\bibitem{Mao97}X.Mao, Stochastic Differential Equations and their Applications, Horwood Publishing,
Chichester, 1997.

\bibitem{Mohammed1984} S.E.A.Mohammed, Stochastic Functional Differential Equations, Research Notes
in Mathematics 99, Pitman Books, London, 1984.

\bibitem{Mohammed1996}S.E.A.Mohammed, Stochastic differential systems with memory.
Theory, examples and applications. In: L. Decreusefond, J. Gjerde,
B. $\O$ksendal, A.S. $\ddot{\mbox{U}}$st$\ddot{\mbox{u}}$nel
(Eds.), Stochastic Analysis and Related Topics VI. The Geilo
Workshop 1996, Progress in Probability,
Birkh$\ddot{\mbox{a}}$user, Basel, 1998, pp.1-77.

\bibitem{SM96}Y.Saito, T.Mitsui, Stability analysis of numerical schemes for
stochastic differential equations, SIAM J. Numer. Anal.
33 (1996) 2254-2267.

\bibitem{Torelli89} L.Torelli, Stability of numerical methods for delay differential equations, J. Comput. Appl. Math.
25 (1989) 15-26.

\bibitem{Wang-Chen-11}W.Wang, Y.Chen, Mean-square stability of semi-implicit Euler
method for nonlinear neutral stochastic delay differential
equations, Appl. Numer. Math. 61 (2011) 696-701.

\bibitem{Wang2008}W.Wang, Numerical analysis of nonlinear neutral functional differential
equations, PhD Thesis, Xiangtan University, 2008. (in Chinese)

\bibitem{WZ2010}W.Wang, C.Zhang, Preserving stability implicit Euler method for nonlinear Volterra and
neutral functional differential equations in Banach space, Numer.
Math. 115 (2010) 451-474.

\bibitem{Wang-Zhang-06}Z.Wang, C.Zhang,  An analysis of stability of Milstein
method for stochastic differential equations with delay, Comput.
Math. Appl. 51 (2006) 1445-1452.

\bibitem{WMS2010}F.Wu, X.Mao, L.Szpruch, Almost sure exponential stability of numerical
solutions for stochastic delay differential equations, Numer.
Math. 115 (2010) 681-697.

\bibitem{ZGH09}H.Zhang, S.Gan, L.Hu, The split-step backward Euler method for
linear stochastic delay differential equations, J Comput. Appl.
Math. 225 (2009) 558-568.


\end{thebibliography}
\end{document}